\documentclass[12pt,a4paper]{article}
\usepackage{amsmath,amssymb,amsthm,graphics,graphicx,epic,eepic,multicol,ascmac}
\usepackage{float}
\usepackage{cite}

\usepackage[all]{xy}

\numberwithin{equation}{section}
\newtheorem{thm}{Theorem}[section]
\newtheorem{prop}[thm]{Proposition}
\newtheorem{lem}[thm]{Lemma}
\newtheorem{rem}[thm]{Remark}
\newtheorem{cor}[thm]{Corollary}
\newtheorem{ex}[thm]{Example}
\newtheorem{con}[thm]{Conjecture}

\theoremstyle{definition}
\newtheorem{defn}[thm]{Definition}

\usepackage{tikz}

\begin{document}

\title{Products of Kirillov-Reshetikhin modules  and maximal green sequences
}

\author{Yuki Kanakubo\thanks{Ibaraki University,
2-1-1, Bunkyo, Mito, Ibaraki, 310-8512,
Japan: {j\_chi\_sen\_you\_ky@eagle.sophia.ac.jp}.}, Gleb Koshevoy\thanks{Institute for Information Transmission Problems Russian Academy of Sciences,
Russian Federation :
{koshevoyga@gmail.com}.} and Toshiki Nakashima\thanks{Division of Mathematics, 
Sophia University, Kioicho 7-1, Chiyoda-ku, Tokyo 102-8554,
Japan: {toshiki@sophia.ac.jp}.} } 
\date{}
\maketitle
\begin{abstract}
We show that a $q$-character of a Kirillov-Reshetikhin module (KR modules) for untwisted quantum affine algebras of simply laced types $A_n^{(1)}$, $D_n^{(1)}$, $E_6^{(1)}$, $E_7^{(1)}$, $E_8^{(1)}$  might be obtained from a specific cluster variable of a seed obtained by applying a maximal green sequence to the initial (infinite) quiver of the Hernandez-Leclerc cluster algebra.  For a collection of KR-modules with nested supports, we show an explicit construction of a cluster seed, which has cluster variables corresponding to the $q$-characters of KR-modules of such a collection. 
We prove that the product of KR-modules of such a collection is a simple module. We also construct cluster seeds with cluster variables corresponding to $q$-characters of KR-modules of some non-nested collections. We make a conjecture that tensor products of KR-modules for such non-nested collections are simple. We show that the cluster Donaldson-Thomas transformations for double Bruhat cells for $ADE$ types can be computed using $q$-characters of KR-modules.
            
\end{abstract}

\section{Introduction
}

Let $\mathfrak g$ be a finite-dimensional simple Lie algebra over $\mathbb C$ and $U_q(\hat{ \mathfrak g})$ be the untwisted quantum affine algebra of $\mathfrak g$. 

Due to Chari and Pressley \cite{CP1}, the isomorphic classes of irreducible finite-dimensional representations of $U_q(\hat {\mathfrak g})$ are parameterized by $n$ -tuples of Drinfeld polynomials, $\mathbf P:=(P_i(u))_{i\in [n]}$, where $n$ is the rank of $ \mathfrak g$, and $[n]=\{1, \cdots , n\}$. The polynomials are of the form
\begin{equation}\label{Dr}
    P_i(u)=\prod_{k=1}^{n_i}(1-c_k^{(i)}u).
\end{equation}
For given Drinfeld polynomials $\mathbf{P}$, let $V(\mathbf P)$ be the corresponding irreducible representation. 

For a pair of Drinfeld polynomials $\mathbf P$ and $\mathbf Q$, let $\mathbf {PQ}:=(P_i(u)Q_i(u))_{i\in I}$. Then $V(\mathbf {PQ})$ is a subquotient of $V(\mathbf P)\otimes V(\mathbf Q)$.

A representation $V(\mathbf P)$ is called the $i$th fundamental representation $V_{\omega _i}(c)$ if $P_i(u)=1-cu$, 
and $P_j(u)=1$ for any $j\neq i$.

The $q$-character $\chi_q$ was introduced by Frenkel and Reshetikhin \cite{FR}. It is an injective homomorphism from the Grothendieck ring ${\rm Rep}\, U_q(\hat{ \mathfrak g})$ of the monoidal category of finite-dimensional representations of $U_q(\hat{ \mathfrak g})$ to the Laurent polynomial ring of infinitely many 
variables $Y_{i, c}$, $i\in [n]$,
$c\in \mathbb C$.
\begin{equation}\label{Dr1}
    \chi_q\,:\, {\rm Rep}\, U_q(\hat{ \mathfrak g})\to \mathbb Z[Y^{\pm 1}_{i,c}]
\end{equation}


A monomial in $\mathbb Z[Y^{\pm 1}_{i,c}]$ is {\em dominant}
if it is a monomial in $Y_{i,c}$, $i\in [n]$, $c\in \mathbb C^{*}$.

Suppose $\mathbf P$ is as in (\ref{Dr}). Then $\chi_q(V(\mathbf P))$ contains the highest weight monomial of $V(\mathbf P)$ 
\begin{equation}\label{Dr2}
M=\prod_i\prod_{k=1}^{n_i}Y_{i,c_k^{(i)}}.
\end{equation}
The highest weight monomials are dominant, and they parameterize the irreducible representations of the untwisted quantum affine algebra of $\mathfrak g$ \cite{FR}. 

Frenkel and Mukhin \cite{FM}
invented an algorithm (FM-algorithm for short)  such that for a given dominant monomial $M$,  it produces a polynomial $\chi_q(M)$. For a special module $V(M)$, they proved that
the result of the algorithm is equal to the $q$-character $\chi_q(V(M))$.
Recall that a module is {\em special}
if its $q$-character contains only one dominant monomial.
The Kirillov-Reshetikhin modules are special \cite{Nakajima, Hernandez1}.

The FM-algorithm is based on the explicit formula of the $q$-characters of the irreducible representations of $U_q(\widehat {\mathfrak{sl}_2})$ (\cite{CP, FR}). 

Another algorithm is due to Nakajima \cite{Nakajima}. For simply laced underlying algebras, Nakajima defined a $t$-analogue of $q$-character as the generating function of the Poincare polynomial of graded quiver varieties. His algorithm is based on the analogue of the Kazhdan-Lusztig polynomials.

Denote by $W_r(c)$ the irreducible representation of $\widehat{{\mathfrak{sl}}_2}$ with the dominant monomial $M=\prod_{k=1}^rY_{cq^{r-2k+1}}$, where $Y_c$ denotes $Y_{1,c}$.
The module $W_r(c)$ is the Kirillov-Reshetikhin module for $\widehat {\mathfrak sl}_2$.

The $q$-character of $W_r(c)$ is 
\begin{equation}\label{sl2-1}
\chi_q(W_r(c))=M(\sum_{i=0}^r\prod_{j=1}^iA_{cq^{r-2j+1}}^{-1}), \quad \mbox{ where }A_b:=Y_{bq^{-1}}Y_{bq}.
\end{equation}
Note that for $c'\neq c$ we get the same result by replacing $Y_c$ by $Y_{c'}$ in (\ref{sl2-1}). 


From \cite{CP1} we know that the $q$-character of any irreducible representation of $\widehat{{\mathfrak{sl}}_2}$ is given by a certain product of $q$-characters of KR modules. Namely, denote by $\Delta_{c, r}$ the set of variables indices in the dominant monomial of $W_r(c)$,
\[
cq^{-r+1}, cq^{-r+3}, \cdots, cq^{r-1},
\]
 and call this set a {\em q-string}. Then
$q$-strings $\Delta_{c, r}$ and
$\Delta_{c', r'}$ are said to be in a {\em general position} if either
$\Delta_{c,r}\cup \Delta_{c',r'}$ is not a $q$-string or one $q$-string is a subset of another. 
Then a given dominant monomial $M$ can be uniquely factorized (up to permutations) into a product of dominant monomials of KR modules with $q$- strings being pairwise in general position. 

Then the $q$-character of $V(M)$ is equal to  the product of $q$-characters of KR-modules for such a set of $q$-strings \cite{CP1}.





Theorem 2.5 in \cite{HL1}, states that the Hernandez-Leclerc category $\mathcal C_l$ is a monoidal category, and its Grothendieck ring $K_0(\mathcal C_l)$ has the structure of a cluster algebra of the finite type $A_l$. The cluster variables of $K_0(\mathcal C_l)$ are the classes
of evaluation modules containing in $\mathcal C_l$. The cluster monomials are equal to the classes of simple modules. Two cluster variables are {\em compatible} (that is they belong to the same cluster seed) if and only if the corresponding strings are in general position.


For ADE-types (we consider simply-laced cases in order to avoid technicalities),  we are interested in the question which collections of KR-modules may correspond to cluster variables of a single seed of the Hernandez-Leclerc cluster algebra. 

In this paper,
first, we obtain a generalization of the above cluster implementation of the Chari and Presley result. 
Namely, we define an analogue of the $q$-string of a KR-module,
a {\em support of its q-character}. We define supports to be in a {\em general position} if they are nested or they are at the distance at least $h/2+1$, where $h$ is the Coxeter number.
We show in Theorems \ref{KRseed} and \ref{main} that, for a collection of KR modules with supports in pairwise general position, there exists a seed whose cluster variables contain the $q$-characters of such a collection (after a monomial transformation of variables),
and the tensor products of the modules of such a collection are simple. 
Note that our Theorem \ref{KRseed}  is related to Theorem 8.1 (a) in \cite{KashiwaraKorea}. 

To prove Theorem \ref{KRstring}, we construct an explicit sequence of mutations in the Hernandez-Leclerc cluster algebra \cite{HL} from the initial seed
to the desired seed such that there is a cluster monomial for this seed, which 
becomes the $q$-character of the tensor products of the modules in a collection after the monomial transformation.

 Second, in Theorem \ref{nonested} we show that there exists a seed of the Hernandez-Leclerc cluster algebra that contains all cluster variables corresponding to a family of $q$-characters of KR-modules whose supports are not in the general position. 

 Third, we use the method of computation of $q$-characters using the source-sink maximal green sequence to get a relation between the cluster Donaldson-Thomas transformation for cluster algebras with the initial seed of the form of the triangular product of a Dynkin quiver with alternating orientation of edges and $A_m$ quiver with the edge directed to the unique sink vertex, and $q$-characters of KR-modules, Theorem \ref{DT-0}. This allows us to use the Frenkel-Mukhin  algorithm or Nakajima algorithm, or a new algorithm proposed in the paper to compute the cluster Donaldson-Thomas transformation. These algorithms compute cluster DT-transformation much faster than cluster computations, since they do not use divisions of polynomials in many variables.

The paper is organized as follows. In Section 2 we remind the reader of the basics of cluster algebra and the triangular product of quivers. In Section 3 we recall the framing of quivers and $c$-vectors and $g$-vectors, and the separation formula. We compute all cluster variables of the $A_{n+1}$ cluster algebra (it is a finite-type cluster algebra) in Proposition \ref{gvectors(n)} and state the cluster interpretation of the result on $q$ strings in general position due to Chari and Pressley. At the end of the section, we define the notion of a maximal green sequence for infinite quivers. In Section 4, we define the source-sink sequence for the triangular product of the Dynkin quiver of ADE types and $A_n$ and prove the level property. By the combination of this level property and Proposition \ref{gvectors(n)}, we compute $g$-vectors of cluster variables in the seed of the cluster algebra with the initial seed being the triangular product, Proposition \ref{slice1}.  
In Section 5, we provide an explicit relation between the cluster variables after applying a reddening sequence to an ambient quiver and its subquiver. In Section 6 we recall necessary facts on the Hernandez-Leclerc cluster algebras.
Section 7 contains the main results, Theorems 7.4, 7.5, 7.9, 7.10, and Conjecture 7.11. In Section 8 we explain how to compute the cluster Donaldson-Thomas transformation for cluster algebras with initial quiver being a triangular product of Dynkin and $A_m$ by using $q$-characters of Kirillov-Reshetikhin modules. %




\section{Cluster algebra}\label{cluster}

\subsection{Quivers}

A quiver $Q=(V_Q,E_Q)$ is a directed graph in the set of vertex $V_Q$ and the set of edges $E_Q\subset V_Q\times V_Q$. As usual $(u,v)\in E_Q$ means that $u$ 
is the head of the edge $e$ and $v$
is the tail.

An ice quiver is a quiver $Q$ such that some non-empty set of vertices $V_Q^{mt}\subset V_Q$ is
distinguished as a set of mutable vertices, 
and the complement $V_Q^{fr}:=V_Q\setminus V_Q^{mt}$ as frozen vertices.
We assign to each vertex $v\in V_Q$, a formal variable $x_v$, such that $(x_v)_{v\in V}$ forms a transcendental basis of the field $\mathbb C((x_v)_{v\in V})$.

\subsection{Cluster seeds and mutations }
Let us briefly recall the cluster algebras (for more details, see
\cite{FZ4}) needed here.  For a positive integer $N$, an ${N}$-regular
tree, denoted by $\mathbb T_N$, whose edges are labeled by $1, \ldots , {N}$, so that the ${N}$ edges emanate
from each vertex with different labels. 
Then an edge of $\mathbb T_N$ is denoted by $t { \rightarrow}_k\, t'$,
indicating that the vertices
$t, t'\in \mathbb T_N$ form an edge $(t,t')$ of $\mathbb T_N$ and $k\in [{N}]$ is the label of this edge.

An $X$-cluster seed is a pair $\Sigma=(\mathbf x, Q)$, where $Q=(V_Q, E_Q)$ is an ice quiver
and $\mathbf x = (x_{j}, \, j\in V_Q)$ is a tuple of variables,
such that the collection $\{x_j| j\in V_Q\}$ generates a field $\mathbf C(x_j,\,j\in V_Q)$.

Let us assign a cluster seed to a root $t_0$ of the tree $\mathbb T_N$ with $N=|V_Q^{mt}|$, we denote this seed
by $\Sigma_{t_0}:=(\mathbf x_{t_0}, Q_{t_0})$. 
An 
X-cluster seed pattern is an assignment of a cluster seed $\Sigma_t=(\mathbf x_t=(x_{j;t})_{j\in V(Q_t)} , Q_t)$ to every vertex $t\in \mathbb T_N$, such
that the seeds assigned to the endpoints of any edge $t\to _v t'$
are obtained from each
other by the seed mutation $\mu_v$, $v\in V_Q^{mt}$. The quiver
mutation $\mu_v$ transforms $Q_t$ into a new quiver $Q_{t'} = \mu_v(Q_t)$ via a sequence of three steps. 

Firstly, for each oriented two-arrow path $u\to v\to w$, $u$, $w\in V_t:=V_{Q_t}$, add a new arrow
$u\to w$. Secondly, reverse the direction of all arrows incident to the vertex $v$. Finally,
repeatedly remove oriented 2-cycles until unable to do so. $E_{Q_{t'}}$ denotes the set of edges of $Q_{t'} = \mu_v(Q_t)$.

For $ X$-variables, we apply 
the cluster mutation rule. Namely, $\mathbf x_{t'}=\mu_v(\mathbf x_t)$, where
$\mu_v(x_{j;t})=x_{j;t}$ if $j\neq v$, and
\begin{equation}\label{mut1}
\mu_v(x_{v;t})=\frac{\prod_{(u\to v)\in E(Q_t) }x_{u;t}+\prod_{(v\to w)\in E(Q_t) }x_{w;t}}{x_{v;t}}.
\end{equation}

A cluster algebra is generated by the $X$-variables of all seeds corresponding to the nodes of $\mathbb T_N$.


A 
$Y$-cluster seed is a pair $\hat \Sigma=(\mathbf y, Q)$,  
where a tuple of variables $\mathbf y = (y_{j}, \, j\in V_Q)$ is
such that the collection $\{y_j| j\in V_Q\}$ generates a field $\mathbf C(y_j,\, j\in V_Q)$.
Y-cluster seeds assigned to the nodes $t\rightarrow_v t'$ are related by a mutation $\hat\mu_v:\hat\Sigma_t\to \hat\Sigma_{t'}$. The quivers are mutated as above. 
For $Y$-cluster variables, we can also apply cluster mutation rule. Namely, $\mathbf y_{t'}=\hat\mu_v(\mathbf y_t)$
\begin{equation}\label{x-cluster}
\hat\mu_v( y_{j,t})=\begin{cases}
\frac 1{y_{j,t}} & if \,j=v\cr
y_{j,t}(1+y_{v,t})^{\# (v,j)}  & \,\,   if \,\,(v,j)\in E\cr

y_{j,t}(1+y_{v,t}^{-1})^{-\# (j,v)} & \,\, if \, \, (j,v)\in E\cr
\end{cases}
\end{equation}
where $\# (i,j)$ denotes the number of edges $(i,j)\in E$.


There is so-called $p$-map which relates $X$-cluster and $Y$-cluster variables
\begin{equation}\label{y-x-trans}
    y_v=\frac{\prod_{v\to u}x_u}{\prod_{w\to v}x_w}.
\end{equation}

\begin{rem}\label{invedges}
Note that for the initial seed $\Sigma_{t_0}^{op}:=(\mathbf x_{t_0}, Q^{op}_{t_0})$ with $Q^{op}_{t_0}=(V_{t_0}, E^{op}_{t_0})$ obtained by reversion of all edges of the initial seed  $\Sigma_{t_0}:=(\mathbf x_{t_0}, Q_{t_0})$,
we get the same 
cluster algebra from $\Sigma_{t_0}^{op}$ as $\Sigma_{t_0}$.
\end{rem}

\subsection{Triangular products
}\label{tri-prod-sec}

The triangle product $Q\boxtimes  R$  of quivers  $Q=(V_Q, E_Q)$ and $R=(V_R, E_R) $ is obtained from their 
direct product $Q \times R$ by adding an arrow from
 $(q,q') \in
V_Q\times V_R$  to $(p,p')\in V_Q\times V_R$  if
$Q$ contains an arrow from $p$ to $q$ and  $R$ contains an arrow from $p'$ to $q'$.

This is a particular case of the definition of a triangular product of valuated quivers in \cite{GK}.
The following is
the triangular product of $A_3^{alt}=v_1\leftarrow v_2\rightarrow v_3$
and $A_3=1\leftarrow 2\rightarrow3$:
\[
\begin{xy}
(0,0) *{\bullet}="3,1",
(-7,5) *{(v_3,1)}="3,1a",
(20,0)*{\bullet}="3,2",
(15,5)*{(v_3,2)}="3,2a",
(40,0)*{\bullet}="3,3",
(47,5)*{(v_3,3)}="3,3a",
(0,-20) *{\bullet}="2,1",
(-7,-15) *{(v_2,1)}="2,1a",
(20,-20)*{\bullet}="2,2",
(26,-25)*{(v_2,2)}="2,2a",
(40,-20)*{\bullet}="2,3",
(47,-15)*{(v_2,3)}="2,3a",
(0,-40) *{\bullet}="1,1",
(-7,-45) *{(v_1,1)}="1,1a",
(20,-40)*{\bullet}="1,2",
(15,-45)*{(v_1,2)}="1,2a",
(40,-40)*{\bullet}="1,3",
(47,-45)*{(v_1,3)}="1,3a",
\ar@{<-} "3,1";"3,2"^{}
\ar@{<-} "3,2";"3,3"^{}
\ar@{<-} "2,1";"2,2"^{}
\ar@{<-} "2,2";"2,3"^{}
\ar@{<-} "1,1";"1,2"^{}
\ar@{<-} "1,2";"1,3"^{}
\ar@{<-} "3,1";"2,1"^{}
\ar@{->} "2,1";"1,1"^{}
\ar@{<-} "3,2";"2,2"^{}
\ar@{->} "2,2";"1,2"^{}
\ar@{<-} "3,3";"2,3"^{}
\ar@{->} "2,3";"1,3"^{}
\ar@{->} "3,1";"2,2"^{}
\ar@{<-} "2,2";"1,1"^{}
\ar@{->} "3,2";"2,3"^{}
\ar@{<-} "2,3";"1,2"^{}
\end{xy}
\]

Picture 1. the triangular product $A^{alt}_3\boxtimes A_3$

\bigskip

Let $G$ be a simply connected connected simple algebraic group or rank $n$,
$B,\ B^-\subset G$ its Borel subgroups, $T:=B\cap B^-$ the maximal torus,
$W={\rm Norm}_G(T)/T$ Weyl group, 
$U$, $U^-$ be unipotent radicals
of $B$, $B^-$,
$A=(a_{i,j})$ the Cartan matrix
of $G$ with an index set $[n]=\{1,2,\cdots,n\}$. 
Let $h$ and $h^\vee$ denote
the Coxeter and dual Coxeter numbers, respectively. 
For the simply-laced case, it holds $h=h^\vee$ and 
the list of $h$ for each type is as below.
\begin{table}[h]
  \begin{tabular}{|c|c|c|c|c|c|} \hline
  \text{type} & $A_n$ & $D_n$ & $E_6$ & $E_7$ & $E_8$ \\ \hline
  ${h}$ & $n+1$ & $2n-2$ & $12$ & $18$  & $30$ \\ \hline
  \end{tabular}
\end{table} 

For each of the $ADE$ types, we consider
a corresponding Dynkin quiver $Q=(V_Q, E_Q)$ with alternating orientation edges.

Below are depicted Dynkin quivers $A_7$, $D_8$ and $E_6$, $E_7$ and $E_8$  with numeration of vertices as in Sage and the alternating orientations of edges.
\[
\begin{xy}
(-8,0) *{A_{7} : }="A1",
(0,0) *{\bullet}="1",
(0,-3) *{1}="1a",
(10,0)*{\bullet}="2",
(10,-3)*{2}="2a",
(20,0)*{\bullet}="3",
(20,-3)*{3}="3a",
(30,0)*{\bullet}="4",
(30,-3)*{4}="4a",
(40,0)*{\bullet}="5",
(40,-3)*{5}="5a",
(50,0)*{\bullet}="6",
(50,-3)*{6}="6a",
(60,0)*{\bullet}="7",
(60,-3)*{7}="7a",
\ar@{<-} "1";"2"^{}
\ar@{->} "2";"3"^{}
\ar@{<-} "3";"4"^{}
\ar@{->} "4";"5"^{}
\ar@{<-} "5";"6"^{}
\ar@{->} "6";"7"^{}
\end{xy}
\]
\[
\begin{xy}
(-8,0) *{D_{8} : }="A1",
(0,0) *{\bullet}="1",
(0,-3) *{1}="1a",
(10,0)*{\bullet}="2",
(10,-3)*{2}="2a",
(20,0)*{\bullet}="3",
(20,-3)*{3}="3a",
(30,0)*{\bullet}="4",
(30,-3)*{4}="4a",
(40,0)*{\bullet}="5",
(40,-3)*{5}="5a",
(50,0)*{\bullet}="6",
(50,-3)*{6}="6a",
(60,8)*{\bullet}="7",
(60,5)*{7}="7a",
(60,-8)*{\bullet}="8",
(60,-11)*{8}="8a",
\ar@{<-} "1";"2"^{}
\ar@{->} "2";"3"^{}
\ar@{<-} "3";"4"^{}
\ar@{->} "4";"5"^{}
\ar@{<-} "5";"6"^{}
\ar@{->} "6";"7"^{}
\ar@{->} "6";"8"^{}
\end{xy}
\]
\[
\begin{xy}
(-8,0) *{E_{6} : }="A1",
(0,0) *{\bullet}="1",
(0,-3) *{1}="1a",
(10,0)*{\bullet}="2",
(10,-3)*{3}="2a",
(20,0)*{\bullet}="3",
(20,-3)*{4}="3a",
(20,10)*{\bullet}="2uu",
(20,13)*{2}="2u",
(30,0)*{\bullet}="4",
(30,-3)*{5}="4a",
(40,0)*{\bullet}="5",
(40,-3)*{6}="5a",
\ar@{<-} "1";"2"^{}
\ar@{->} "2";"3"^{}
\ar@{<-} "3";"4"^{}
\ar@{->} "4";"5"^{}
\ar@{<-} "3";"2u"^{}
\end{xy}
\]
\[
\begin{xy}
(-8,0) *{E_{7} : }="A1",
(0,0) *{\bullet}="1",
(0,-3) *{1}="1a",
(10,0)*{\bullet}="2",
(10,-3)*{3}="2a",
(20,0)*{\bullet}="3",
(20,-3)*{4}="3a",
(20,10)*{\bullet}="2uu",
(20,13)*{2}="2u",
(30,0)*{\bullet}="4",
(30,-3)*{5}="4a",
(40,0)*{\bullet}="5",
(40,-3)*{6}="5a",
(50,0)*{\bullet}="6",
(50,-3)*{7}="6a",
\ar@{<-} "1";"2"^{}
\ar@{->} "2";"3"^{}
\ar@{<-} "3";"4"^{}
\ar@{->} "4";"5"^{}
\ar@{<-} "3";"2u"^{}
\ar@{<-} "5";"6"^{}
\end{xy}
\]
\[
\begin{xy}
(-8,0) *{E_{8} : }="A1",
(0,0) *{\bullet}="1",
(0,-3) *{1}="1a",
(10,0)*{\bullet}="2",
(10,-3)*{3}="2a",
(20,0)*{\bullet}="3",
(20,-3)*{4}="3a",
(20,10)*{\bullet}="2uu",
(20,13)*{2}="2u",
(30,0)*{\bullet}="4",
(30,-3)*{5}="4a",
(40,0)*{\bullet}="5",
(40,-3)*{6}="5a",
(50,0)*{\bullet}="6",
(50,-3)*{7}="6a",
(60,0)*{\bullet}="7",
(60,-3)*{8}="7a",
\ar@{<-} "1";"2"^{}
\ar@{->} "2";"3"^{}
\ar@{<-} "3";"4"^{}
\ar@{->} "4";"5"^{}
\ar@{<-} "3";"2u"^{}
\ar@{<-} "5";"6"^{}
\ar@{->} "6";"7"^{}
\end{xy}
\]
Picture  2.\label{pic1}
\medskip
 
 
 

For Dynkin quiver $Q$ with alternating orientation of edges, let $V^{sc}=\{p_1, \cdots, p_k\}$ be the set of source vertices and $V^{sk}=\{q_1, \cdots, q_l\}$ be the set of sink vertices so that $V_Q=V^{sc}\cup V^{sk}$.
Each edge takes the form $(p_i, q_j)$ for some $i\in [k]$ and $j\in [l]$. Note that for any pair of vertices $p_i$, $p_{i'}\in V^{sc}$, we have zero entries in the Cartan matrix, $a_{p_i, p_{i'}}=0$ and the same $a_{q_j, q_{j'}}=0$ for any pair of sink vertices in $V^{sk}$.

For the alternating edge orientation of Dynkin quiver $Q$ of $DE$ types and $A_{2k+1}$,  the 
corresponding reduced decomposition of $w_0$ 
takes the form \[
w_0=C^{h/2},
\]
where $C$ is a reduced decomposition of the Coxeter element
\begin{equation}\label{cox}
C=q_1\cdots q_l\cdot p_1\cdots p_k.
\end{equation} 
For $A_{2k}$, we take $w_0=C^{\frac{h-1}2}q_1\cdots q_l$.
We define a Coxeter word $C^*:=q_1^*\cdots q_l^*\cdot p_1^*\cdots p_k^*$, where $i^*\in V_Q$ is defined by
$\alpha_{i^*}=-w_0\alpha_i$.


Following Hernandez and Leclerc \cite{HL}, we are interested in cluster algebras with the infinite quiver
\begin{equation}\label{HLquiver}
Q\boxtimes A_\infty,
\end{equation}
where $A_\infty$ is the semi-infinite line quiver with vertices labeled by positive integers and edge orientations such that there are no source vertices and  the unique sink vertex is labeled by $1$.

The quiver (\ref{HLquiver}) is almost the same quiver as $G^-$ in \cite{HL}.
  The latter quiver is the triangular product $Q\boxtimes A_{-\infty}$, where $A_{-\infty}$ is a quiver with nodes labeled by $-1, -2, \cdots$, such that
  the unique sink vertex $-1$ and no source vertices (see Section \ref{HL}). 

For $Q=A_{2n+1}, DE$ types, one can regard the triangular product 
$Q\boxtimes A_\infty$ as quiver for the infinite word $(C^{h/2}(C^*)^{h/2})^\infty$ in the sense of \cite{BFZ3}, namely, this word has the property that any subword of consecutive letters of length $l(w_0)$ is a reduced decomposition of $w_0$ (subsection 4.3 of \cite{KashiwaraKorea}).  
For $Q=A_{2n}$, one has to consider the word $(C^{h/2})^\infty$, for which any word of length $l(w_0)$ is a reduced decomposition.
 



\section{Framing, green sequences, reddening sequences}\label{basicfram}

\begin{defn}
	The \emph{framed quiver $\overline{Q}$} associated to a quiver $Q$ is obtained by adding for each vertex $v\in V_Q$ a frozen vertex $\bar{v}$ and an arrow $v\rightarrow \bar {v}$. We set $\overline{V}:=\{\overline{v}|v\in V_Q\}$.
	Dually, we define the \emph{co-framed quiver~$\underline{Q}$} by adding arrows $\overline {v} \rightarrow v$ instead.
\end{defn}

For a quiver $Q=(V_Q,E_Q)$, we set $\overline Q=(V_{\overline{Q}}:=V_Q\cup \overline V, E_{\overline{Q}}:=E_Q\cup \{(v,\overline v), v\in V)$ and $\underline Q=(V_Q\cup \overline V, E_Q\cup \{(\overline v, v), v\in V)$. 
Note that, for each frozen vertex $v\in V_Q^{fr}$, we also add a new frozen vertex\footnote{One can regard the frozen vertices $V_Q^{fr}$ as mutable, however we are nether going to mutate at them.}. The frozen vertices of $\overline Q$ are defined as $V_{\overline{Q}}^{fr}:=
V_Q^{fr}\cup \overline{V}$.

The {\em framed} seed $\overline\Sigma=(\mathbf x, \overline {\mathbf x}, \overline Q)$ is associated to the seed $\Sigma$, where $\overline {\mathbf x}=( x_{\bar v})_{\bar v\in \bar V}$. 
\begin{defn} For an initial framed seed $\overline \Sigma$,    a corresponding  tree $\mathbb T_{N}$,  and a vertex $t$ of the tree,
	a mutable vertex $v\in \overline Q_t$ is called \emph{green} if no arrow departing from a frozen vertex $w\in\overline{V}$ targets $v$. Similarly, a mutable vertex $v\in \overline Q_t$ is called \emph{red} if no arrow from $v$ to a frozen vertex $w\in\overline {V}$.
\end{defn}

    
Let $\mathcal{S}=(v_{i_1}, \ldots , v_{i_m})\in  (V^{mt}_{Q_{t_0}})\times (V^{mt}_{Q_{t_1}}) \times \cdots \times (V^{mt}_{Q_{t_{m-1}}})$ be a sequence of mutable vertices, where $Q_{t_i}$ is a quiver associated with the vertex $t_i$ of $\mathbb T_N$ 
obtained by the first $i$ mutations following $\mathcal S$. Any mutable vertex in \[
\overline {Q}\mathcal{S}:=\mu_{v_{i_m}}\cdots \mu_{v_{i_1}}(\overline Q)
\] has precisely one of the two colors, red or green \cite{DWZ}.

Following Keller \cite{Keller}, we define certain sequences relating~$\overline{Q}$ and $\underline {Q}$. 
\begin{defn}
	A sequence $\mathcal{S}=(v_{i_1}, \ldots , v_{i_m})$ is called \emph{green} sequence if the vertices $v_{i_j}\in \overline{Q} \mathcal{S}_{j-1}$ are green for all $j<n$ ($\mathcal S_k\subset \mathcal S$ is constituted of the first $k$ vertices). A sequence $\mathcal{S}$ is called \emph{reddening } if all mutable vertices are red in $\overline{Q}\mathcal{S}$. A green reddening sequence is called a \emph{maximal green sequence} (MGS for brevity). 
\end{defn}
By \cite[Proposition 2.10]{BDP} for each reddening sequence $\mathcal{S}\subset Q$ there exists a unique isomorphism of ice quivers $$\sigma_{\mathcal{S}}:( \overline{Q}\mathcal{S})  \xrightarrow{\sim} \underline{Q}$$
 fixing the frozen vertices and sending a non frozen vertex $v$ to $\sigma (v)$, $v\in V^{mt}_{\overline Q}$.

\begin{rem}
 Let, for an initial seed $\Sigma=(\mathbf x, R)$, there are two reddening sequences $\mathcal S_1$ and $\mathcal S_2$, and let $\Sigma_1=(\mathbf x ^1=(x_v^1)_{v\in V_R}, R\mathcal S_1)$ and $\Sigma_2=(\mathbf x^2=(x_v^2)_{v\in V_R}, R\mathcal S_2)$ be seeds after applying mutations by these sequences.
 Then, for $v\in V_R^{mt}$, we have 
 \begin{equation}\label{twoseqid}
   x^1_v=x^2_{\sigma_{\mathcal S_1}(
 \sigma_{\mathcal S_2}^{-1}(v))}.  
 \end{equation}
This follows because the $g$-vectors of $x_v^1=-e_{\sigma_{\mathcal S_1}v}$ and
$g$-vectors of $x_v^2=-e_{\sigma_{\mathcal S_2}v}$,
and $g$-vectors determine cluster variables \cite{CKQ}.
 
\end{rem}
\subsection{$c$-vectors and BPS charges}
BPS charges, $c$-vectors and $g$-vectors are defined with respect to an initial quiver.
For an initial  quiver $Q$ and a vertex $v\in V_Q^{mt}$, the edge $(v,\bar v)$ of $\overline Q$ defines the initial $c$-{\em vector} of $v$.
Namely, we identify the initial $c$-vector of a vertex $v$ with a basic unit vector
in $\mathbb Z^{ V^{mt}}$, \[(e_v)_{ w}=\begin{cases}
1\mbox{ if } w=v, \cr 
0 \, \mbox{ if } w\in V^{mt}\setminus \{v\}. \cr
\end{cases}\] 
For a
sequence $\mathcal S=(v_{i_1}, \ldots, v_{i_m})$, 
let us pick the subsequence ${\mathcal S}_{j-1}$,
then $c$-vector of a vertex $v$  of the quiver 
$\overline Q{\mathcal S}_{j-1}$ are defined as follows:
Consider the set of frozen vertices
$C(v)\subset \overline V$ to which $v$ is joined in the quiver $\overline Q{\mathcal S}_{j-1}$  (with multiplicities). 
Then the $c$-vector of  $v$ is
 $c_{v}=\sum_{\overline w\in C(v)} e_w$ if 
 $v$ is green and
 $c_{v}=-\sum_{\overline w\in C(v)} e_w$ if $v$ is red.
Then the mutation at the vertex $v=v_{i_{j}}$  of the quiver $\overline Q\mathcal S_{j-1}$ changes $c$-vectors by the following rule 
 \begin{equation}\label{cvectors1}
\mu_v(c_w)=\begin{cases}
-c_v\mbox{ if }w=v\cr
c_w+\sum_{(w, v)\in E_{j-1}}c_{v} \mbox{ if $v$ is green}\cr
c_w+\sum_{(v,w)\in E_{j-1}} c_{v} \mbox{ if $v$ is red}\cr
\end{cases}
\end{equation}
where $E_{j-1}$ denotes the set of edges of the quiver $Q\mathcal S_{j-1}$.

For a $c$-vector 
$\sum_{(v,j), \, v\in C, j\in J}c_{(v,j)}e_{(v,j)}$ with
$c_{(v,j)}\in \mathbb Z_{\neq 0}$, $C\subset V_Q$, $J\subset \mathbb Z_+$, 
the set $\{(v,j)| v\in C, j\in J\}$
is the {\em support set} of the $c$-vector.
Following \cite{Xie}, we define the BPS charges.
\begin{defn}
For a green sequence of mutations $\mathcal S=\{v_{i_1}, \cdots, v_{i_m}\}$ and each $j=1, \cdots, m$,
the $c$-vector at the vertex $v_{i_j}$ of the
quiver $\overline{Q}\mathcal S_{j-1}$ is called a BPS-charge. The sequence of such $c$-vectors forms BPS-charges for $\mathcal S$.
\end{defn}

Of interest are BPS-charges for maximal green sequences, since they are finite sequences of $c$-vectors.

\begin{defn}
A sequence $\mathcal S=\{v_{i_1}, \cdots, v_{i_m}\}$ is a {\em  sink sequence}, if each $v_{i_j}$ is a green sink vertex in $ Q\mathcal S_{j-1}$, $j=1, \ldots, m$.
A sequence $\mathcal S=\{v_{i_1}, \cdots, v_{i_m}\}$ is a {\em  source sequence}, if each $v_{i_j}$ is a green source vertex in $Q\mathcal S_{j-1}$, $j=1, \ldots, m$.
\end{defn}

For any Dynkin quiver $Q$, there exist sink and source maximal green sequences.


\begin{prop}\label{source1}
For any finite acyclic quiver $R$ (no oriented cycles), there exists a source maximal green sequence.\end{prop}
{\em Proof}.
 In fact, mutation at a green source vertex $v$ changes the colour only at this vertex and all other vertices preserve their colours. 
The remaining subgraph of green vertices is acyclic as well, and therefore there is a source vertex of this subgraph. Since we mutate at source vertices,  there are no edges directed from red vertices to green ones, hence we get that such a  source
vertex is a  green source vertex of the whole graph. We mutate at 
such a green source vertex and so on, and end up with a  source maximal green source.\hfill $\Box$

\begin{ex}
For $A_5$ of Picture 2, we have the following source maximal green sequences:
\[ 
\{2,4,1,3,5\},\, \{4,2,1,3,5\}, \{2,4,3,1, 5\},\, \{4,2,3,1, 5\}, \{2,4,3,5,1\},\, \{4,2,3,5,1\},\]
\[
\{2,4,1,5,3\},\, \{4,2,1,5,3\}, 
\{2,4,5,1,3\},\, \{4,2,5, 1,3\}, \{2,4,5,3,1\},\, \{4,2,5, 3,1\},\]
\[\{2,1, 4,3,5\},\, \{4,5, 2,3, 1\},
\{2,1, 4,5,3\},\, \{4,5, 2,1, 3\}.
\]
\end{ex}

Concerning sink maximal green sequences, we have
\begin{prop}\label{sink1}
Let $Q$ be a Dynkin quiver with some edge orientation. Then there exists a sink maximal green sequence.
\end{prop}
{\em Proof}. It is easy to check that for an edge oriented Dynkin quiver $Q$, the sink sequence which gives a reduced decomposition of $w_0$ related to the edge orientation $Q$ is a sink maximal green sequence.
\hfill $\Box$

Note that for a given orientation, two such decompositions are related by 2-moves only.
\medskip

\begin{ex}\label{anex}
One can check that, for the Dynkin quiver $A_n$ with usual edges oriented such that $1$ is the single sink vertex,  the following sequences
are sink maximal green sequences indeed.
The first one is
\[
\mathcal S_{A_n}=\{1, 2, 1, 3, 2, 1, 4, 3, 2, 1, \ldots, n, n-1, n-2, \ldots, 1\}.\]
The second one is
\[
\mathcal S^{HL}_n =\{1, 2, \ldots, n, 1,  2, \ldots, n-1, 1, 2, \ldots, n-2, \ldots, 1, 2, 1\}.\]
\end{ex}


\begin{lem}\label{BPS1}
The BPS-charges for the sequence $\mathcal S_{A_n}$ are

\[e_1, e_1+e_2, e_2, e_1+e_2+e_3, e_2+e_3, e_3, \ldots, e_1+e_2+e_3+\cdots +e_n, \]
\begin{equation}\label{order1}
e_2+e_3+\cdots +e_n, e_3+\cdots +e_n, e_{n-1} +e_n, e_n.
\end{equation}
\end{lem}
{\em Proof}. 
For a Dynkin quiver $Q$, the set of $c$-vectors for sink maximal green sequence is in  bijection with the set  positive roots of the corresponding Lie algebra \cite{BDP}. Namely, for a sink maximal green sequence $\mathcal S=\{v_{i_1}, \cdots, v_{i_m}\}$,
the order of the roots corresponds to 
 the reduced decompositions of the longest element $w_0=s_{i_1}\cdots s_{i_m}$. For the sequences $\mathcal S_{A_n}$, the decomposition is
\begin{equation}\label{i0}
 s_1s_2s_1s_3s_2s_1s_4s_3s_2s_1\cdots s_ns_{n-1}\cdots s_1.
\end{equation}
This claim is verified by induction on $n$. In fact, for the sequence $\mathcal S_{A_n}$, by induction, we have validity of (\ref{order1}) for $n-1$.
Before the last round of the mutations we have the following $c$-vectors: at the vertex $j\in [n-1]$, $c$-vector is $-e_{n-1-j+1}$ (the inverse permutation) and these vertices are red, and at the vertex $n$ is $e_1+\cdots+e_n$.

Applying the mutation at $n$, the $c$ vector of $n$ is $- (e_1+\cdots+e_n)$ and the vertex $n-1$ becomes green with the $c$ vector
\[
e_2+\cdots +e_n,\]
since there is a cyclic triangle $(n-1)\to n\to \overline 1\to (n-1)$ before mutation at $n$ therefore, after mutation at $n$, the vertex $(n-1)$ will join with all frozen of $n$ except $1$. 
Hence, after mutation at $(n-1)$ it will have the $c$-vector $- (e_2+\cdots+e_n)$, the vertices $(n-2)$ and $n$ will turn green with the $c$-vectors $(e_3+\cdots+e_n)$ and $e_1$, respectively. Continuing, we get the claim.
 \hfill $\Box$


\subsection{$g$-vectors and separation formula}

Denote by $C=C^{B,t_0}_t$ the $c$-matrix of column $c$-vectors, $c_v(t)$, $v\in V_{mt}$.

Define a matrix $G_t$, the $g$-matrix, from the equation
\begin{equation}\label{gvectorsdef}
G^T_t=(C^{(B^T,t_0)}_t)^{-1}
\end{equation}
where $B^T$ denotes transpose matrix. 
Then columns of the matrix $G_t$, $g_v(t)$ are said to be $g$-{\em vectors}.
This way of defining $g$-vectors is not standard but is correct because of the main results in \cite{NZ} and the sign coherence of $c$-vectors.

The separation formula allows us to express a cluster variable $x_{v,t}$ using its $g$-vectors $g_v(t)$ and $F$-polynomials: It takes the form
\begin{equation}\label{separ1}
x_{v,t}={\mathbf x}_{t_0}^{g_v(t)}F_{v,t}(\mathbf y_{t_0}).
\end{equation}

The basic role of $g$-vectors is that they
provide natural labels of cluster variables, namely $g$-vectors characterize cluster variables $x_{v,t}$, $t\neq t_0$ \cite{FZ4}.

In relations to $q$-characters of KR-modules, we are interested in computing $g$-vectors for nodes  of the quiver $A_{n+1}\mathcal S_{A_n}$, $n=1, \ldots$.

\begin{prop}\label{gvectors(n)}
Let $\Sigma_0=(x_1,\cdots,x_n,x_{n+1};A_{n+1})$ be an initial seed with the quiver of type $A_{n+1}$ defined in Example \ref{anex}
and let $x^{(n+1)}_{a}$ $(a=1, \ldots, n)$ denote the cluster variable of the seed $\Sigma_0 \mathcal S_{A_n}=
(x_1^{(n+1)}, \cdots, x_n^{(n+1)}, x_{n+1}; A_{n+1}\mathcal S_{A_n})$.
Then $x^{(n+1)}_{n+1-a}$
 is equal to the following Laurent polynomial 
\begin{equation}\label{variable-an}
x^{(n+1)}_{n+1-a}=\frac{x_{n+1}x_{n}\cdots x_{a+1}(1+ \frac{x_{a-1}}{x_{a+1}}(1+\frac{x_{a}}{x_{a+2}}(1+ \frac{x_{a+1}}{x_{a+3}}(\cdots+\frac{x_{n-2}}{x_{n}} (1+\frac{x_{n-1}}{x_{n+1}})))))}{x_nx_{n-1}\cdots x_a}.
\end{equation}
In particular, its $g$-vector is $g( x^{(n+1)}_{n+1-a})=e_{n+1}-e_a$ and its BPS charge is $e_a+\cdots +e_n$.
\end{prop}
{\em Proof}. 
We  verify (\ref{variable-an}) by induction.

Namely, (\ref{variable-an})  takes the form 
\begin{equation}\label{Fform1}
\frac{x_{n+1}}{x_{a}}(1+ y_a+y_ay_{a+1}+\cdots +y_ay_{a+1}\cdots y_n)
\end{equation}
since \[
y_j=\frac{x_{j-1}}{x_{j+1}}, \, j=1, \ldots, n.
\]

Observe that it holds
\begin{equation}\label{y-x}
y_a\cdots y_n=\frac{x_{a-1}}{x_n}\frac{x_a}{x_{n+1}}.
\end{equation}
We have to compute the variables of the seed $\Sigma'_{(n+2)}:=\mu_1\mu_2\cdots \mu_{n+1}\tilde \Sigma_{n+2}$, where 
the seed $\tilde \Sigma_{n+2}=(x^{(n+1)}_1,\cdots, x^{(n+1)}_n, x_{n+1}, x_{n+2}; A_{n+2}^{n+1})$,
 and  $A_{n+2}^{n+1}$ being the Dynkin quiver $A_{n+2}$ has the edge orientation such that the vertex $n+1$ is the single sink vertex.

Then 
\[
x^{(n+2)}_{n+1}=\mu_{n+1}(x_{n+1})=
\frac{ x_{n+2}\frac{x_{n+1}}{x_{1}}(1+ y_1+y_1y_{2}+\cdots +y_1y_{2}\cdots y_n)+1}{x_{n+1}}=
\]
\begin{equation}\label{eq3}
\frac{x_{n+2}}{x_{1}}(1+ y_1+y_1y_{2}+\cdots +y_1y_{2}\cdots y_n)+\frac 1{x_{n+1}}=
\end{equation}
\begin{equation}\label{eq4}
\frac{x_{n+2}}{x_{1}}(1+ y_1+y_1y_{2}+\cdots +y_1y_{2}\cdots y_n + y_1y_{2}\cdots y_{n+1}). 
\end{equation}
(\ref{eq4}) is obtained by substituting $y_1y_{2}\cdots y_{n
+1}=\frac 1{x_{n+1}}\frac{x_1}{x_{n+2}}$ due to (\ref{y-x}).

Because in the quiver $A_{n+2}$, we should label nodes by $n+2-a=n+1-(a-1)$, and 
the quiver $\mu_{n+2-a}\cdots \mu_{n+2-1}A_{n+2}$ has a unique sink vertex at $n+2-(a+1)$,
we have 
\[
x^{(n+2)}_{n+2-(a+1)}
=\frac{
x^{(n+1)}_{n-a}x^{(n+2)}_{n+2-a}+1}
{x^{(n+1)}_{n+1-a}}
\]
\[=\frac{
x^{(n+1)}_{n-a}\frac
{x_{n+2}}{x_{a}}
(1+ y_{a}+y_{a}y_{a+1}+\cdots +y_{a}y_{a+1}\cdots y_{n}+y_{a}y_{a+1}\cdots y_{n+1})  +1
}
{
x^{(n+1)}_{n-a+1}
}.
\]
Then we get 
\[
x^{(n+2)}_{n+2-(a+1)}= \frac
{x_{n+2}}{x_{a+1}}
(1+ y_{a+1}+y_{a+1}y_{a+2}+\cdots +y_{a+1}y_{a+2}\cdots y_{n})+
\]
\begin{equation}\label{eq5}
\frac
{x_{n+2}}{x_{a+1}}\frac{\frac{x_{a+1}}{x_{n+2} }\frac{x_{a}}{x_{n+1}} +(1+y_{a+1}+y_{a+1}y_{a+2}+\cdots +y_{a+1}\cdots y_n)y_a(y_{a+1}\cdots y_{n+1})}{1+y_{a}+y_{a}y_{a+1}+\cdots +y_{a}\cdots y_n}.
\end{equation}
Because of (\ref{y-x}), we get that the last line of (\ref{eq5}) is equal to \[
\frac {x_{n+2}}{x_{a+1}}y_{a+1}y_{a+2}\cdots y_{n}y_{n+1}.\]
Thus, we get (\ref{variable-an}).

For computing $g$-vectors, we use (\ref{Fform1}) and the separation formula (\ref{separ1}). 
The BPS charge at the vertex $n+1-a$ is computed in Lemma \ref{BPS1}.\hfill $\Box$ 
\begin{rem}\label{sl-2}
The $q$-characters for Kirillov-Reshetikhin $U_q(\hat sl_2)$-modules (\cite{CP}) one can  obtain from (\ref{Fform1}) by the following change of variables
\[
Y_{q^k}:=\frac{x_k}{x_{k-1}}.
\]
Then \[
y_a=\frac 1{Y_{q^a}Y_{q^{a+1}}}=:A_a^{-1}
\]
and (\ref{Fform1}) takes the form 
\begin{equation}\label{Fform2}
x^{(n+1)}_{n+1-a}=Y_{q^{n+1}}\cdots Y_{q^{a+1}}(1+ A_a^{-1}+A_a^{-1}A_{a+1}^{-1}+\cdots +A_a^{-1}A_{a+1}^{-1}\cdots A_n^{-1})
\end{equation}
that is noting else as $q$-character of KR-module $W_{n-a+1}({q^{n+a+2}})$
(see (\ref{sl2-1})).

\begin{prop}\label{strings2} Let $\mathcal L;=\{W_r(q^s)$ $|$ $r\in A$, $s\in S\}$ be a collection of $\widehat{ sl_2}$ KR-modules with $q$-strings in a general position. Then, for some $n$ there exists a sequence of mutations $\mathcal S$ of $A_{n+1}$-cluster algebra with frozen vertex $n+1$, such that
for some set $J$, the cluster variables $x^{(n+1)}_j$ of the seed $\mu_\mathcal S\Sigma_0$, $j\in J $, become $q$-characters of  modules of $\mathcal L$ 
 after the above monomial change of variables.
\end{prop}
We get this proposition as a particular case $Q=A_1$ of Theorem \ref{KRseed}.
\end{rem}

\subsection{Maximal green sequences for infinite quivers}

We define a maximal green/reddening sequence for an infinite quiver as injective limit.
\begin{defn}\label{mgsinfty}
For an infinite quiver $Q$, a sequence  $\mathcal S$
is called a maximal green sequence if there exists an increasing sequence of subquivers $Q^1\subset Q^2\subset \cdots$ such that $\cup Q^j=Q$ and a collection 
of mutations $\mathcal S^i$ of $Q^i$, $i=1, \ldots$, such that
\[
\mathcal S^j\mbox{ is a maximal green/reddening sequence for } Q^j.
\]
and the sequence $\mathcal S^{j-1}$ is the beginning part of the sequence $\mathcal S^j$, and we have $\mathcal S=(\mathcal S^1, \mathcal S^2\setminus \mathcal S^1, \cdots
\mathcal S^{i+1}\setminus \mathcal S^i \ldots)$.
\end{defn}

For an infinite quiver $Q$ that has a maximal green sequence, we define the BPS charge by the same rule as for the finite case.
Here is an important example. 

\begin{ex}
Consider the quiver $A_\infty$ with the vertex set $\{1, 2, 3, \ldots\}$ and the arrows $(i+1, i)$, $i=1,2, \ldots$.

Below we depict 
$\overline A_\infty$.

\begin{xy} 0;<1pt,0pt>:<0pt,-1pt>:: 
(0,83) *+{1} ="0",
(60,83) *+{2} ="1",
(120,83) *+{3} ="2",
(180,83) *+{4} ="3",
(240,83) *+{5} ="4",
(300,83) *+{\cdots} ="4a",
(2,0) *+{1'} ="5",
(60,2) *+{2'} ="6",
(120,3) *+{3'} ="7",
(180,5) *+{4'} ="8",
(240,6) *+{5'} ="9",
(300,6) *+{\cdots} ="9a",
"1", {\ar"0"},
"0", {\ar"5"},
"2", {\ar"1"},
"1", {\ar"6"},
"3", {\ar"2"},
"2", {\ar"7"},
"4", {\ar"3"},
"3", {\ar"8"},
"4", {\ar"9"},
"4a", {\ar"4"},
\end{xy}
\bigskip

We consider the following sink sequence of mutations
\[
\mathcal S_A:=1, 2,1, 3, 2, 1, 4,3,2,1, \cdots.\]

\begin{lem}\label{infA1}
The sequence $\mathcal S_A$ is a maximal green for the infinite quiver $A_\infty$.
\end{lem}
{\em Proof}. The sequence $1,2,1,3,2,1,\cdots , n, n-1, \cdots, 2, 1$ is a sink MGS for
each $n\in\mathbb{Z}_+$. Then, 
by the sequence of quivers $A_n$ and the mutation sequences
$\mathcal S_{A_n}$, we get a maximal green sequence for $A_\infty$.
\hfill $\Box$

\end{ex}

For an infinite quiver $Q$ that has a maximal green sequence $\mathcal S$,  the corresponding BPS sequence of $c$ vectors and $g$ vectors are gradually defined for each $Q^i$ and $\mathcal S^i$, $i=1, \ldots$.



\section{Source-sink maximal green sequences for triangular products}
\subsection{Source-sink sequences and the level property}

For the triangular product $Q\boxtimes D$, where $Q$ is an acyclic quiver  and $D$ is a Dynkin quiver,   a source-sink  mutation sequence $Sc(Q)\times Sk(D)$ was defined in \cite{GK}  as follows: 
let $Sc(Q)=v_{1}, \cdots, v_{n}$ be a source maximal green sequence for $Q$ and $Sk(D)= i_1, \ldots , i_N$ be  a sink  maximal green sequence  for $D$.  
Then the source-sink sequence $Sc(Q)\times Sk(D)$ is 
\[(v_{1},i_1),  (v_{2},i_1),\cdots, (v_{n},i_1),  
\]
\[(v_{1},i_2),  (v_{2},i_2),\cdots, (v_{n},i_2),  
\]
\[\cdots\]
\[(v_{1},i_N),  (v_{2},i_N),\cdots, (v_{n},i_N). 
\]

Theorem 3.8 in \cite{GK} states that the source-sink sequence  $Sc(Q)\times Sk(D)$   is a
maximal green sequences  for the triangular product $Q\boxtimes D$.

For example for $A_2\boxtimes A_3$, $A_2=v_1\leftarrow v_2$, $A_3=1\leftarrow 2\leftarrow 3$,
the source-sink sequence $Sc(A_2)\times Sk(A_3)$  with $Sc(A_2)=2,1$ and $Sk(A_3)=1,2,1,3,2,1$ is the sequence
\[
(v_2,1), (v_1, 1), (v_2, 2), (v_1, 2), (v_2, 1), (v_1, 1),  (v_2,3), (v_1,3), (v_2, 2), (v_1, 2), (v_2, 1), (v_1, 1).
\]

\subsection{Level property}
We consider the {\em level property} in detail for the triangular product $Q\boxtimes A_n$ with a Dynkin quiver of type $ADE$ and alternating edge orientation.
Let $p_1, \cdots, p_k$ and $q_1, \cdots, q_l$ denote the source vertices and the sink vertices $Q$, respectively.

Of our interest is the following source maximal green sequence for $Q$
\begin{equation}\label{sourceQ1}
Sc(Q)=p_1, \cdots, p_k, q_1, \cdots, q_l.\end{equation}

Then, for a sink maximal green sequence $Sk(A_n)=\mathcal S_{A_n}= i_1, \ldots , i_N$, we consider the subsequence $\mathcal S_j$ of $S_C(Q)\times S_{A_n}$ in the previous subsection:

\[
\mathcal S_{j}:=(p_{1},i_1),  (p_{2},i_1),\cdots, (p_{k},i_1),  (q_{1},i_1),  (q_{2},i_1),\cdots, (q_{l},i_1),
\]
\[(p_{1},i_2),  (p_{2},i_2),\cdots, (p_{k},i_2), (q_{1},i_2),  (q_{2},i_2),\cdots, (q_{l},i_2),  
\]
\[\cdots\]
\begin{equation}\label{seqtr1} (p_{1},i_{j-1}),  (p_{2},i_{j-1}),\cdots, (p_{k},i_{j-1}),  (q_{1},i_{j-1}),  (q_{2},i_{j-1}),\cdots, (q_{l},i_{j-1}). 
\end{equation}
For example, when $n=5$, $Q=A_4$, applying mutations $\mathcal S_3(v_1,4),(v_3,4),(v_2,4),(v_4,4)$, we get the following quiver:
\begin{figure}[H]
\centering
  \includegraphics[width=8cm, height=6cm]
  {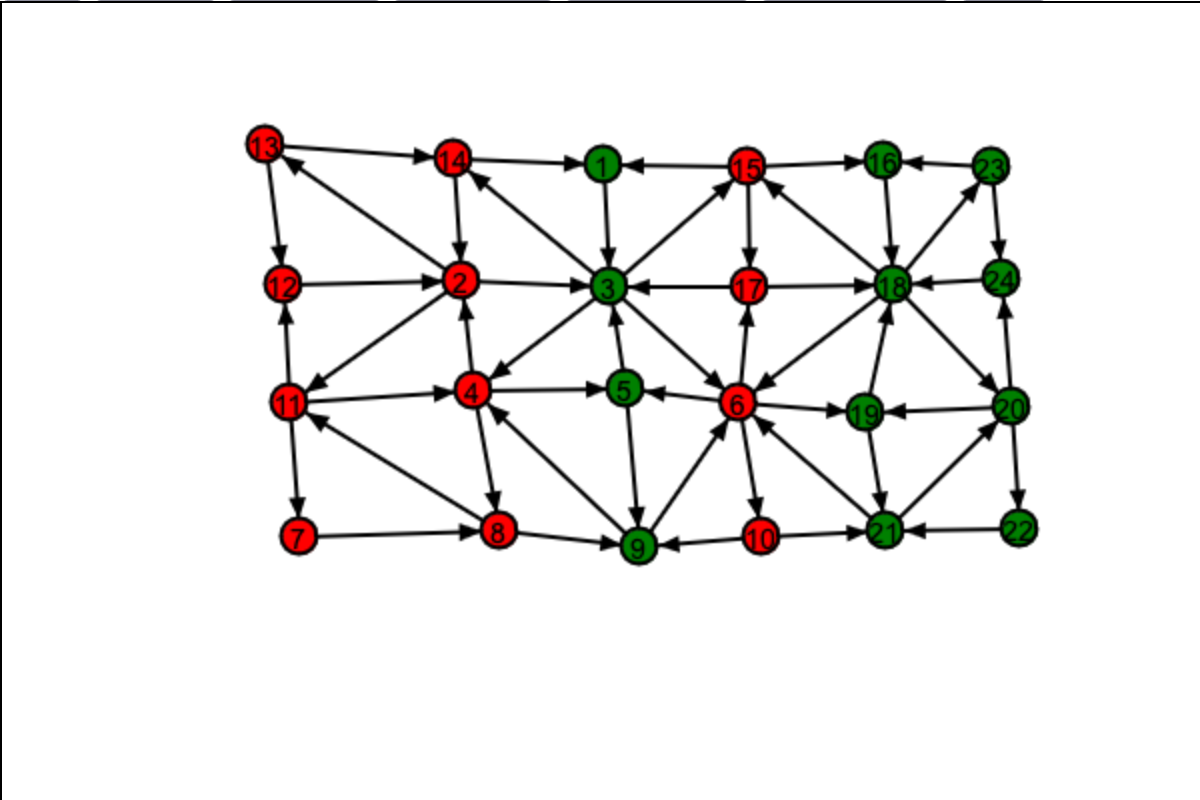}
  \end{figure}  
 $\mbox{ Picture 3\label{c-vect-pic}}$
 
\hspace{-8mm}
Note that in our sequence we will mutate the vertex $(v_1,3)$ next time, which is green and has only two ingoing edges $((v_1,2),(v_1,3))$
and $((v_1,4),(v_1,3))$.

Because of the definition of the triangular product, we get the following:  following any source-sink sequence, when we mutate at a green vertex of the form $(\cdot, i_j)$, such a vertex has only two ingoing edges
\begin{equation}\label{ingoing}
((\cdot, i_j-1), (\cdot, i_j)) \mbox{ and } ((\cdot, i_j+1), (\cdot, i_j)).
\end{equation}
Note that a similar argument appears in the proof of
Proposition 3.6 of \cite{GK}.

\begin{defn}    
For the quiver $Q\boxtimes A_m$, a sequence of mutations has the {\em level property}, if $c$-vectors at each mutable vertex $(v, j)$, $v\in V^{mt}_Q$, $j\in [m]$, have the support
of the set $\{(v, i)| i\in [m]\}$.
\end{defn}
We show above that the source-sink sequence has the level property.
\begin{lem}\label{levelsub}
   For a quiver $Q\boxtimes A_m$ and $m'<m$, 
  consider, for the subquiver $Q\boxtimes (A_m\setminus A_{m'})$, the mutation sequence  $Sc(Q)\times \mathcal S_{A_m\setminus A_{m'}}$.
  Then this sequence has the level property.
  \end{lem}
  {\em Proof}. 
In case of $m'=0$, 
let $(v_r,i_j)$ be a vertex in $Sc(Q)\times S_{A_m}$, where
$v_r\in\{p_1,\cdots,p_k,q_1,\cdots,q_l\}$.
Since $Sc(Q)\times S_{A_m}$ is a source-sink maximal green sequence,
just before mutating at $(v_r,i_j)$,
it holds
\begin{itemize}
\item[(1)] $(v_r,i_j)$ is green.
\item[(2)] $(v_r,i_j)$ is a source in subquiver 
$(v_1,i_j), (v_2,i_j), ..., (v_n,i_j)$.
\item[(3)] By (\ref{ingoing}), the vertex has only two ingoing edges
\[
((v_r, i_j-1), (v_r, i_j)) \mbox{ and } ((v_r, i_j+1), (v_r, i_j)).
\]
\end{itemize}
By these, we see that $Sc(Q)\times S_{A_m}$ has the level property.
For general $m'$,
the claim follows from the level property of source-sink sequences with $m'=0$, since the shown above properties of mutation in green source vertices remain valid for $m'>0$. 
\hfill $\Box$  
\begin{prop}\label{slice1}
For a quiver $Q\boxtimes A_{m+1}$, let $\Sigma_0=(x_{(v_i, j)}, \, i\in [n], \, j\in [m+1]\}, Q\boxtimes A_{m+1})$ be the initial seed.
Let $x^{(m+1)}_{(v_i,a)}$, $i\in [n]$, $a\in [m]$ and $x_{(v_i,m+1)}$, $i\in [n]$ denote the cluster variables in the seed $\Sigma_0 (Sc(Q)\times \mathcal S_{A_m})$. 
Then 
\begin{equation}\label{sep-product1}
x^{(m+1)}_{(v_i,m+1-a)}= \frac{x_{(v_i,m+1)}}{x_{(v_i,a)}}F(\mathbf y)
\end{equation}
where $F(\mathbf y)$ is $F$-polynomial corresponding to $x^{(m+1)}_{(v_i,m+1-a)}$.
\end{prop}

{\em Proof}. 
Due to the level property of the maximal green source-sink sequence, we can reduce the calculation of the $c$ vectors for the sequence $Sc(Q)\times \mathcal S_{A_m}$ to $\mathcal S_{A_{m}}$  for each $v_i$, $i\in [n]$.

For $\mathcal S_{A_m}=v_{i_1}, \cdots, v_{i_N}$, the $c$-matrix for $Sc(Q)\times  S_{A_m}$ is a block matrix and each block is a c-matrix for the mutation sequences $(v_i,i_1),\cdots,(v_i,i_N)$ in the fiber quiver $\{v_i\}\times A_{m+1}$, for each $i\in [n]$.

Then by ({\ref{gvectorsdef}), the $g$-matrix of $(Q\boxtimes A_{m+1})Sc(Q)\times \mathcal S_{A_m}$is a block matrix.  
Hence, we get the statement from Proposition
\ref{gvectors(n)}.
\hfill $\Box$

\begin{rem}\label{propslice}
In Proposition \ref{slice1}, the quiver $Q\boxtimes A_{m+1}$ is the ice quiver with frozen vertices $\{(v_i, m+1))\}$, $i\in [n]$. For a quiver $Q\boxtimes A_{L}$ with $L\ge m+1$, one may regard vertices
$\{(v_i, m+1))\}$, $\cdots$, $\{(v_i, L))\}$, $i\in [n]$, to be frozen. Let $\Sigma ^L$ be the corresponding seed. Then  (\ref{sep-product1}) holds for $\Sigma^L Sc(Q)\times \mathcal S_{A_m}$.
\end{rem}

Of our interest are the triangular products of the form $Q\boxtimes A_{\infty}$, where $Q$ is a Dynkin quiver of type $ADE$ with alternating orientation of edges, and $A_\infty$ is Dynkin quiver of $A$-type with infinite set of nodes labeled by positive integers $\mathbb N=\{1, 2, \ldots \}$ and unique sink vertex $1$.

The point is that Hernandez and Leclerc in \cite{HL} have considered exactly such quivers for applying cluster algebra algorithm to compute $q$-characters of KR-modules (see Section \ref{HL}).


\subsection{Source-sink for $Q\boxtimes A_\infty$}

For $A_\infty$,  following Lemma \ref{infA1}, we consider the maximal green sink sequence
\[
\mathcal S_A:=\lim_n \mathcal S_{A_n}=1, 2,1, 3, 2, 1, \cdots.\]

\begin{prop}\label{MGSbasic}
For the infinite quiver
$
Q\boxtimes A_\infty $, 
the source-sink sequence  \[Sc(Q)\times \mathcal S_A,\]
   is a maximal green sequence. For each $m$, Proposition \ref{slice1} holds true with the replacement  the quiver $Q\boxtimes A_{m+1}$ by the quiver $Q\boxtimes A_{\infty}$.
\end{prop}
{\em Proof}. 
 Follows Proposition \ref{slice1} and Remark \ref{propslice}. 
\hfill $\Box$

\section{Quivers and subquivers}

We consider a quiver $Q\boxtimes A_K$ and, for $k<K$, we consider the subquiver $Q\boxtimes A_{[k,K]}$, where $A_{[k,K]}$ is a subquiver of $A_K$ with the vertices $k, k+1, \cdots, K$.

 Consider two sequences of mutations, $\mathcal S_{A_{K-1}}$ and the sequence \[\mathcal S_{A_{[k,K-1]}}=k, k+1, k, k+2, k+1, k, \cdots , K-1, K-2, \cdots, k.
 \]

For the initial seed $\Sigma_0=(\mathbf x, Q\boxtimes A_K)$,
consider seeds $\Sigma_1:=\Sigma_0(Sc(Q)\times \mathcal S_{A_K})=
(\mathbf x', (Q\boxtimes A_K)(Sc(Q)\times \mathcal S_{A_K}))$,
and $\Sigma_2:=\Sigma_0(Sc(Q)\times \mathcal S_{A_{[k,K]}})= (\mathbf x'', (Q\boxtimes A_K)(Sc(Q)\times \mathcal S_{A_{[k,K]}}))$. The latter seed has the same cluster variables as $\Sigma_0$ has at vertices $(v, j)$, $j<k$. 
We are interested in relations between cluster variables of the seeds
$\Sigma_1$ and $\Sigma_2$.

In the beginning, we consider $\widehat {sl_2}$, that is, $Q$ is the singleton quiver
$A_1$. 
Then we get the following claim:
\begin{lem}\label{mgsmM}
    \begin{equation}
      x'_i=x''_{k+i-1}, \quad i=1, \cdots, K-k.
    \end{equation}
\end{lem}
{\em Proof}. Let us consider another maximal green sequence and use (\ref{twoseqid}). That is, we consider the source sequence $\mathcal K=K-1, K-2, \cdots , 1$. Denote by $\Sigma_3=\Sigma_0\mathcal K=
(\mathbf x^1, A_K\mathcal K)$ the corresponding $X$-cluster seed.
For such a sequence, the isomorphism $\sigma_\mathcal K$ is the identical permutation of mutable vertices. 

The cutting of such a sequence,  $\mathcal K'=K-1, K-2, \cdots , k$, is a maximal green sequence for the subquiver with vertices $k, k+1, \ldots, K-1$ (recall that the vertex labeled by $K$ is frozen). 
Denote by $\Sigma_4=\Sigma_0\mathcal K'=
(\mathbf x^2, A_K\mathcal K')$ the corresponding $X$-cluster seed.
The isomorphism $\sigma_{\mathcal K'}$ is identical on the vertices $\{k, k+1, \ldots, K-1\}$. 
Then, for $j=K, K-1, \ldots , k$, it holds $x^1_j=x^2_j$.

For the sequences $\mathcal S_{A_{K-1}}$ and $\mathcal S_{A_{[k,K-1]}}$, the isomorphisms are
$j\mapsto K-j$ ($j=1,2,\cdots,K-1$) and $j\mapsto j$ ($j=1,\cdots,k-1$), $j\mapsto K+k-j-1$ ($j=k,k+1,\cdots,K-1$), respectively.
By (\ref{twoseqid}), we obtain
$x'_j=x^1_{K-j}$ $(j=1,\cdots,K-1)$ and $x''_j=x^2_{K+k-j-1}$ $(j=k,k+1,\cdots,K-1)$.
Therefore, one gets
$x_i'=x^1_{K-i}=x^2_{K-i}=x''_{k+i-1}$.


\hfill $\Box$

A conceptual explanation for Lemma \ref{mgsmM} is as follows. The sink sequence $\mathcal S_{A_{[k,K]}}$ is obtained from $\mathcal S_{A_{K}}$ by
omitting vertices $1, \cdots, k-1$ of the latter one. 
Note that the former sequence is a maximal green sequence for the subquiver $A_{[k, K]}$.

Consider a general case, let $R$ be a quiver with a reddening sequence $\mathcal S$, let $R'$ be a subquiver such that a sequence $\mathcal S |_{R'}$  obtained from $\mathcal{S}$ by omitting vertices outside $R'$ is a reddening sequence for $R'$, let $\sigma $ be the automorphism of vertices $R$ for the reddening sequence $\mathcal S$ and $\sigma' $ be that for $R'$ and $\mathcal S|_{Q'}$.

Then it holds
\begin{prop}\label{Fpolysub}
 For a vertex $v\in V_{R'}$, the F-polynomial for $x'_{\sigma'^{-1}(v)}=\mu_{\mathcal S|_{R'}}(x_{\sigma'^{-1}(v)})$ is obtained form the F-polynomial for
 $x_{\sigma'^{-1}(v)}=\mu_{\mathcal S}(x_{\sigma^{-1}(v)})$ by substituting $y_w=0$ for $w\in V_R\setminus V_{R'}$.
\end{prop} 


{\em Proof}. The proof follows the same idea of degeneration of scattering diagrams as in the proof of Proposition 4.3 in \cite{Wang}. 

 Let $\mathfrak D_R$ be a scattering diagram for the quiver $R$, and we choose function $1+ty_w$ instead of $1+y_w$ for $w$th initial wall with $w\notin V(R')$. Let $p$ be a path that goes from the positive chamber to the negative chamber of $\mathfrak D_R$.

 Let $v\in V_Q$ and $\sigma (v)\in Q'$. Then $x'_{\sigma ^{-1}v}$ is equal to the outcome of the wall-crossing transformation of $x_w$ along the path $\gamma$. This gives F-polynomial in the separation formula for the variable $x'_{\sigma ^{-1}v}$. 
 
 Then we set  $t=0$, by the same arguments as in \cite{Wang} 
 we get the F-polynomial in the separation formula for
 $x'_{(\sigma')^{-1}v}$ by sending $y$-variables to zero outside $R'$, $y_w=0$, $w\in V(R)\setminus V(R')$. 
 \hfill $\Box$
 
\begin{rem}\label{sl2sub}  
Consider the case in Lemma \ref{mgsmM}.
The inverse permutation
$\sigma (i)=K-i$ is the automorphism for $\mathcal S_{A_{K-1}}$ and
$\sigma'(k+i-1)=K-i$, $i=1, \ldots , K-k$, $\sigma'(i)=i$, $i<k$ is the
automorphism for $\mathcal S_{A_{k,K-1]}}$. We regard the latter automorphism as the automorphism of the subquiver $A_{[k,K-1]}$.
Then, by Proposition \ref{Fpolysub}, the F-polynomial of $x'_{k+i-1}$ is obtained from F-polynomial of $x_{i}$, $i=1, \cdots, K-k$, by substituting $y_w=0$ for vertices outside $[k,K]$. It is easy to see that for $Q=A_1$, there are no such $y$-variables. This is exactly stated in Lemma \ref{mgsmM}.
\end{rem}

From the level property, Lemma \ref{levelsub},   Proposition \ref{Fpolysub} and Remark \ref{sl2sub}, we get the following.

\begin{prop}\label{xsub}
For $Q$ of an $ADE$-type and the seeds $\Sigma_1$ and $\Sigma_2$, any $v\in V_Q$ and $1\le i<K-k-h/2$, it holds
\[
x''_{v, k+i-1}=x'_{v, i}
\]
\end{prop}
{\em Proof}. The F-polynomial in (\ref{sep-product1}) could only have $y$-variables labeled by $(v, j)$ with $a-h/2\le j\le m$ (in our case one has to set $m=K-1$ and $a=K-i$). This follows from the relation of such polynomials and $q$-characters of KR-modules (see Proposition \ref{KRcluster} below).
The $q$-character for KR-module $V(M)$ is a partial sum of the product of $q$-characters for fundamental modules with dominant monomials running the set of variables in $M$. F-polynomial of the $q$-character for fundamental module with the dominant monomial $\frac{x_{v,l+1}}{x_{v,l}}$ could have $y$-variables with indices of the set $\{(v,j)| v\in V_Q, l-h/2\le j\le l\}$. (For $AD$-cases this follows from \cite{Nakajima}, and for exceptional cases, one can verify this by computations
using, for example, the Frenkel-Mikhin algorithm.)   
Therefore, by 
\[
\frac{x_{v,K}}{x_{v,K-i}}=\prod_{l=K-i}^{K-1}
\frac{x_{v,l+1}}{x_{v,l}},
\]
the F-polynomial for $x'_{v, i}$ has $y$-variables whose labels
are in $\{(v, j)| K-i-\frac{h}{2}\le j\le K-1 \}$.
In particular, for
$i\le K-k-h/2$, the labels are in $\{(v, j)| k\le j\le K-1 \}$. Hence, for $i\le K-k-h/2$,
the F-polynomial for
$x''_{v, k+i-1}$ has the same set of $y$-variables in the F-polynomial for
$x'_{v, i}$. In fact, for $i<K-k$, we have $\sigma (i)=\sigma' (k+i-1)$. Therefore, for such $i$'s, by Proposition \ref{Fpolysub} we know that the $F$-polynomial for $x''_{v, k+i-1}$
is obtained from $F$-polynomial for $x'_{v,i}$ by substituting $y_{v', j}=0$, with
$v'\in V_Q$, $j<k$. 

For $i<K-k-h/2$, we show above that the set of such $y$-variables is empty.
This implies the proposition.
\hfill $\Box$

\section{Hernandez-Leclerc cluster algebra and the cluster algorithm for $q$-characters of KR-modules}\label{HL}
\subsection{Monomial changes of variables}

For quantum affine algebra of type $A_n^{(1)}$, $D_n^{(1)}$, $E_6^{(1)}$, $E_7^{(1)}$, $E_8^{(1)}$,  we delete the vertex labeled by zero, and  consider
the corresponding $ADE$-type Dynkin quiver $Q=(V_Q, E_Q)$ with alternating orientation of  edges and let $V^{sc}=\{p_1, \cdots, p_l\}$ be the subset of source vertices and $V^{sk}=\{q_1, \cdots, q_m\}$ be the subset of sink vertices,
$V_Q=V^{sc}\cup V^{sk}$.
We keep three sets $V_Q=\{v_1, \cdots, v_n\}$ and the set of sources $V^{sc}=\{p_1, \cdots, p_l\}$ and sinks $V^{sk}=\{q_1, \cdots, q_m\}$, $l+m=n$.
Following Hernandez and Leclerc \cite{HL} we considered a 
cluster algebra $\mathcal A_Q$ with the initial seed $\Sigma_0=(\mathbf x, Q\boxtimes A_\infty )$.  And let $\hat\Sigma_0=(\mathbf y, Q\boxtimes A_\infty )$
be the corresponding  initial $Y$-cluster seed.

 The cluster algebra $\mathcal A_Q$ is the $\mathbb Q$-subalgebra of the field of rational functions $\mathbb Q(\mathbf x)$ generated by all the elements obtained from the initial seed $\Sigma_0$ by a finite sequence of seed mutations.

Following \cite{HL}, we consider a change of coordinates for computing $q$-characters from cluster variables 
\begin{equation}\label{tx} t_{v_j, 1}:=x_{v_j, 1}, \quad t_{v_j, k}:=\frac{x_{v_j, k}}  {x_{v_j, k-1}}, 
\, j\in [n],\, k=2, 3\ldots, 
\end{equation}
Then, the $y$-variables of the seed $\hat\Sigma_0$ take the following form: For $i\in[l]$ and $j\in[m]$, it follows from
(\ref{y-x-trans}) that
\[
y_{p_i, k}:=\frac{x_{p_i,k-1}\prod _{(p_i,q_j)\in E(Q)}x_{q_j,k}}
{x_{p_i,k+1}\prod _{(p_i,q_j)\in E(Q)}x_{q_j,k-1}
},
 \, k=2, \cdots\]
 \[
y_{p_i, 1}:=\frac{\prod _{(p_i,q_j)\in E(Q)}x_{q_j,1}}
{x_{p_i,2}}
\]
 \[
y_{q_j, k}:=\frac{x_{q_j,k-1}\prod _{(p_i,q_j)\in E(Q)}x_{p_i,k+1}}
{x_{q_j,k+1}\prod _{(p_i,q_j)\in E(Q)}x_{p_i,k}
},
 \, k=2, \cdots\]
 \[
y_{q_j, 1}:=\frac{\prod _{(p_i,q_j)\in E(Q)}x_{p_j,2}}
{x_{p_i,2}\prod _{(p_i,q_j)\in E(Q)}x_{p_j,1} }.
\]

Rewriting $y$-variables in $t$-variables gives us 
\[
y_{p_i, k}:=t_{p_i, k}^{-1}\cdot t_{p_i, k+1}^{-1}\prod_{q_j\,:\, (p_i, q_j)\in E(Q)}t_{q_j,k}, \]
\[ y_{q_j, k}:=t_{q_j, k}^{-1}\cdot t_{q_j, k+1}^{-1}\prod_{p_i\,:\, (p_i, q_j)\in E}t_{p_i,k+1 }.
\]

\begin{rem}
Note that in comparing with notations of variables $Y_{t, q^s}$ in \cite{HL},  for a source vertex $v$,  $t_{v, k}$ corresponds to $Y_{v,q^{-2k+2}}$, and for a sink vertex $w$,  $t_{w, k}$ corresponds to $Y_{w,q^{-2k+1}}$, and $y_{v_i,k}$ are the affine analog of $e^{\alpha_{v_i}}$ operators $A_{v_i,k}$. 

\end{rem}

 \subsection{The Kirillov-Rechetikhin modules}\label{KRHL}
 Let $U_q(\hat{\mathfrak g})$ be the quantum enveloping algebra of 
 $\hat{\mathfrak g}$, which is the untwisted affine Lie algebra of the simple Lie algebra $\mathfrak g$. We consider $A_n^{(1)}$, $D_n^{(1)}$, $E_6^{(1)}$, $E_7^{(1)}$ or $E_8^{(1)}$ type.
 
 Frenkel and Reshetihkin \cite{FR} have attached to every complex finite-dimensional  representation $M$ of 
$U_q(\hat{\mathfrak g})$ 
a $q$-character $\chi_q(M)$.   
Irreducible finite-dimensional representations of
$U_q(\hat{\mathfrak g})$
can be parametrized by the highest dominant monomial of their $q$-character \cite{FR}.
By definition, the $q$ -character
$\chi_q(M)$ is a Laurent polynomial with positive integer coefficients in the infinite set of variables
$Y=\{Y_{i,c}| i\in I, c\in \mathbb C^*\}$, where $I$ denotes the set of simple roots of the underlying Lie algebra $\mathfrak g$.
Following Hernandez and Leclerc \cite{HL}, we are interested only in polynomials involving the variables
$\{Y_{i, q^k}\}$, $i\in I$, $k\in\mathbb Z_-$.

In the notations of the previous subsection, we consider polynomials in variables
\[
t_{v_i, j},\, i\in [n],\ j=1, 2, 3,\ldots .\]

 Thus, $q$-characters will be Laurent polynomials in the variables
\[
\mathbf t= (t_{(v_i, j)},\, i\in [n], j=1, 2, 3\ldots ).
\]
 
Let $M$ be a dominant monomial in the variables of $\mathbf t$
, that is, a monomial with nonnegative exponents. We denote by $V(M)$ the corresponding irreducible representation of
$U_q(\hat{\mathfrak g})$ and by $\chi_q(M) =\chi_q(V(M))$ its $q$-character.

If $M$ is of the form
\[
\prod_{j=1}^{k}t_{v, r+j} 
\]  
where $v\in V_Q$, then $V(M)$ is called a {\em Kirillov-Reshetikhin module}, and we denote\footnote{
Since the vertices of the Hernandez-Leclerc quiver are labeled with negative powers of $q$,  there is a difference in notations here comparing to Section 2.3.2 in \cite{HL}. That is, a KR module $W^{(v)}_{k, r+k}$ is denoted in \cite{HL} by $W_{k,-2r-2k+2}^{(i)}$ when $v$ is a source and
by $W^{(i)}_{k, -2r-2k+1}$ when $v$ is a sink, where $i$ the Dynkin quiver node corresponding to $v\in V_Q$. 
} it by
\[W_{k,r+k}^{(v)},
\]
where $v\in V_Q$, $t_{v, r+k}$ is the rightmost variable of $M$, and $k$ is the length of $M$.
We call the interval $[r+1, r+k]:=\{r+1, \cdots, r+k\}$ the {\em support of the dominant monomial } $M=\prod_{j=1}^{k}t_{v, r+j}$.
 Note that the support of a dominant monomial does not depend on a node of $Q$.
For $k=1$ , $W^{(v)}_{1,r}=V(t_{(v,r)}$ is called a fundamental module.

The cluster algorithm for computing the $q$-characters of KR-modules invented by Hernandez and Leclerc \cite{HL} runs as follows.
Recall that we denoted by $h$ the Coxeter number.

Let us run $d\ge h/2$ times the following infinite mutation sequence
\[
\{(p_1, 1), (p_2, 1), \cdots, (p_k, 1),  (q_1, 1), (q_2, 1), \cdots, (q_l, 1),\]
\[ (p_1, 2), (p_2, 2), \cdots, (p_k, 2), 
 (q_1, 2), (q_2, 2), \cdots, (q_l, 2)\]
 \[,\cdots
\]
\begin{equation}\label{HLmodif} (p_1, s), (p_2, s), \cdots, (p_k, s), 
 (q_1, s), (q_2, s), \cdots, (q_l, s),
\end{equation}
\[\cdots
\]

Denote by $\Sigma_{d\infty}$ the seed obtained from the initial seed $\Sigma_0$ by applying the mutation sequence (\ref{HLmodif}) $d$ times.

\begin{rem}
Due to the commutativity of mutations at the source vertices, that is, the following two sequences of mutations terminate at the same seed,  
\[
(p_1, 1), (p_1, 2),\cdots, (p_1, s), (p_2, 1), (p_2, 2),\cdots, (p_2, s), 
\]
and 
\[
(p_1, 1), (p_2, 1), (p_1, 2), (p_2, 2)\cdots, (p_1, s), (p_2, 1), (p_2, s), s=1, \cdot . 
\]
The above sequence of mutations is equivalent to that in Section 3.1 \cite{HL}.
\end{rem}

Then we can state Theorem 3.1 in \cite{HL} as the following statement.
\begin{thm}
For $d\ge h/2$ and    
for any $k\ge 1$, the cluster variables  at the node $(v_i, k)$  of the seed $\Sigma_{d\infty}$ becomes the $q$-character of KR-module  
\[W_{k,d+k}^{(v_i)}
\]
after  change of $x$-variables to $t$-variables by (\ref{tx}).
\end{thm}

From (\ref{tx}) there holds
\begin{equation}
    \prod_{j=1}^{k}t_{v_s, r+j} =\frac{x_{v_s, r+k}}{x_{v_s, r}}, 
\end{equation}
By Proposition 4.16 in \cite{HL}
$g$-vector of the cluster variable corresponding to $q$-character $W_{k, k+d}^{(v_i)}$, $v_i\in V_Q$, is
\[
e_{v_i, d+k}-e_{v_i, d}.
\]

In other words, for the  Hernandez-Leclerc cluster algebra $
\mathcal A_Q$,  a cluster variable whose 
 $g$-vector is
\begin{equation}\label{HLgvector}
e_{v_i, d+k}-e_{v_i, d}.
\end{equation}
 gives the $q$-character of $W_{k, k+d}^{(v_i)}$, for $d\ge h/2$ and truncated
$q$-character for $d< h/2$, after change of variables ({\ref{tx}}).

From this and since $g$-vectors characterize the cluster variables \cite{CKQ}, we get the following:
\begin{cor}\label{gq}
If the $g$-vector of a cluster variable $x$ of the cluster algebra $\mathcal{ A}_Q$  is equal to $e_{v_i, d+k}-e_{v_i, d}$, $d\ge h/2$. Then after change of variables (\ref{tx}), we get  the $q$-character of KR module $$W_{k, k+d}^{(v_i)}.$$
\end{cor}
\section{Maximal green sequences and $q$-characters of products KR-modules}
\subsection{Nested collections}

Combining (\ref{sep-product1}) in Proposition \ref{slice1} and (\ref{HLgvector}), we get the following:
\begin{prop}\label{KRcluster}
For $b\ge  1$ and $L> b+h/2$, and $a\le b$, the cluster variable
\[
x^{(L)}_{(v_i, a)}
\]
 of the seed \, $\Sigma_0 Sc (Q)\times \mathcal S_{A_{L-1}}$
 gives the $q$-character of $W_{a,L}^{(v_i)}$ after change of variables (\ref{tx}) .
 \end{prop}

 {\em Proof}. By (\ref{sep-product1}), we get that the $g$-vector of $x^{(L)}_{v_i, a}$ is 
 \[
e_{(v_i, L)}- e_{(v_i, L-a)}.
 \]
From this and Corollary \ref{gq} we get the Proposition.
 \hfill $\Box$

  \begin{rem}\label{trunc1}
  Since we consider a semi-infinite quiver $Q\boxtimes A_\infty$, for a KR-module $W_{a,L}^{(v)}$ with  $L\le a+h/2$, we can get $q$-character $\chi_q(W_{a,L}^{(v_i)})$ by shifting variables, and  without shifting we get the truncated $q$-characters. 
  \end{rem}

 \begin{defn}
A collection $\mathbf N$ of intervals of $\mathbb Z$, is {\em nested } if, for any intervals $[a+1,a+b]$ and $[c+1,c+d]\in \mathbf N$, it holds either $[a+1,a+b]\subset [c+1,c+d]$ or $[c+1,c+d]\subset [a+1,a+b]$.
\end{defn}
For a nested collection $\mathbf N$, and for quantum affine algebra $\mathfrak g$ of one of $A_n^{(1)}$, $D_n^{(1)}$, $E_6^{(1)}$, $E_7^{(1)}$, $E_8^{(1)}$ types, let 
\[
\mathcal M_\mathbf N=\{W^{(v)}_{b, a+b}| v\in V_Q, \, [a+1,a+b]\in \mathbf N\}
\]
 be the 
 collection of KR-modules whose supports of dominant monomials are in
$\mathbf N$.

\begin{thm}\label{KRseed} 
  Let $Q$ be the Dynkin quiver of the corresponding $ADE$ type with alternating orientation of edges and let $\Sigma_0=(\mathbf x, Q\boxtimes A_\infty)$ be an initial seed of the Hernandez-Leclerc cluster algebra $\mathcal A_Q$. Let a finite collection $\mathbf N$ and $\mathcal M_\mathbf N$ be as above.
Then there exists a finite sink green sequence of mutations $\mathcal S_\mathbf N$, 
  such that the seed $\Sigma'=\Sigma_0\mathcal S_\mathbf N$ contains cluster variables whose $g$-vectors take the form
  \[
  e_{(v, a+b)}-e_{(v, a)},
  \]
 $[a+1, a+b]\in \mathbf N$, $v\in V_Q$. 
\end{thm}

 Because of Theorem \ref{KRseed}, for a collection $\mathbf N$ such  that all intervals are in the right side of $h/2$, products of $q$-characters of KR-modules of $\mathcal M_\mathbf N$ correspond to cluster monomials, one may ask on simplicity of the tensor products of such KR-modules.  
 

\begin{thm}\label{main}
Under the assumptions of Theorem \ref{KRseed},
the tensor product 
\[
\otimes_{W\in \mathcal M_\mathbf N} W,
\]
is a simple module for any ordering of the modules in the tensor product.
\end{thm}
{\em Proof of Theorem \ref{main}}. For a bipartite Coxeter word  
$C=\prod_jq_j\prod_ip_i$, let $C^*$ be the word defined in the subsection \ref{tri-prod-sec}.
Let
$\hat{\underline w_0}$ be the infinite word
$(C^{h/2}(C^*)^{h/2})^\infty$ or $(C^{h/2})^\infty$ at
the end of subsection \ref{tri-prod-sec}, which has
the property that any subword of consecutive letters of length $l(w_0)$ is a reduced decomposition of $w_0$.
Let us write $\hat{\underline w_0}=(\imath_l)_{l\in\mathbb{Z}}$.
For $k\in\mathbb{Z}$, if $\imath_k$ is a source in $Q$ then we set $\tilde{p}_k=-2\lfloor \frac{-k-1}{n}\rfloor$ and if $\imath_k$ is a sink then set $\tilde{p}_k=-1-2\lfloor \frac{-k-1}{n}\rfloor$,
where $\lfloor c \rfloor$ denotes maximal integer smaller than or equal $c$ for $c\in\mathbb{R}$.
In the notation and terminology of \cite{KashiwaraKorea},
we take an admissible sequence $\mathfrak{s}=(\imath_k,\tilde{p}_k)_{k\in\mathbb{Z}}$.
Note that the partial sequence
$\mathfrak{s}=(\imath_k,\tilde{p}_k)_{k\in\mathbb{Z}_{<0}}$ coincides with the set of vertex of the quiver $G^-$ in \cite{HL}.
Following Theorem 6.14 of \cite{KashiwaraKorea}, one can associate a KR-module
\[
M[a,b]=M^{\mathfrak{s}}[a,b]
\]
to each $i$-box $[a,b]$. 
For an interval $[-b,-a]$ with $-b<-a<0$ and $\imath_{-a}=\imath_{-b}=i$, it holds
$M[-b,-a]=W^{(i)}_{|[-b,-a]|,q^{\tilde{p}_{-b}}}$ by Remark 6.15 of \cite{KashiwaraKorea}, which coincides
with $W^{(i)}_{|[-b,-a]|,\lfloor \frac{b-1}{n}\rfloor+1}$ in our notation. Here, $|[-b,-a]|$ is the number
of elements in $\{s\in\mathbb{Z}|-b\leq s\leq -a,\ \imath_s=\imath_{-a}=\imath_{-b}\}$.


We identify $\{(\imath_k,-\tilde{p}_k)|k\in\mathbb{Z}_{<0}\}$ with $\mathbb{Z}_{<0}$ by $(\imath_k,-\tilde{p}_k)\leftrightarrow k$.
Then, 
the requirement that $[-b, -a]$ is a $i$-box is equivalent to   
\begin{itemize}
\item{(WQ)} \quad $-a$ corresponds to a vertex $(v_i,j)$ and $-b$ corresponds to $(v_i, j+k)$ for the vertex $v_i=\iota_{-a}$, and positive $j$ and non-negative $k$. 
\end{itemize}


Note that in case of (WQ), $M[-b,-a]=W^{(i)}_{|[-b,-a]|,\lfloor \frac{b-1}{n}\rfloor+1}=W^{(i)}_{k+1,j+k}$.
Two $i$-boxes $[-b_1,-a_1]$ and $[-b_2,-a_2]$ commute (Definition 4.18 in \cite{KashiwaraKorea}) if 
the corresponding intervals $[j_1, k_1+ j_1]$ and $[j_2, k_2+ j_2]$ defined in WQ are nested, where $(v_{i_s},j_s)$ and $(v_{i_s},k_s+j_s)$ correspond to $-a_s$ and $-b_s$, $s=1, 2$. 
 Theorem \ref{main} follows from
 Corollary 3.16
 of \cite{KKKO} and
 Theorem 4.21, Theorem 5.5 in \cite{KashiwaraKorea} and Theorem 1.1 in \cite{Hernandez}.
 \hfill $\Box$
 \medskip

Note that, for $ADE$ types, Theorem 7.7 in \cite{OT} stating that tensor product of Kirillov-Reshetikhin modules corresponding to a cluster monomial of a seed obtained by applying the Hernandez-Leclerc sequence of mutations is simply follows from Theorem \ref{main}. In fact $g$-vectors of cluster variables of such a seed take the form 
\[
e_{v, k}-e_{v, h/2}, \quad k>h/2
\]
 and therefore their supports are nested.
 \bigskip
 
 {\em Proof of Theorem \ref{KRseed}}

 Let $\mathbf N=\{[r_i+1, r_i+k_i]\}_{i=1, \ldots m}$ be nested collection and $r_1\le r_2\le \ldots \le r_m$. 
For any $ADE$ Dynkin quiver $Q$, we set $\mathcal S_\mathbf N$ to be source-sink sequence $Sc(Q)\times \mathcal S_{A_\mathbf N}$.
It suffices to establish a required sequence of mutation for an adjoint 
pair of modules, say $W_{k_{m-1},r_{m-1}+k_{m-1}}^{(v)}$ and $W_{k_m,r_m+k_m}^{(v)}$.

 First, we consider the case 
 $Q$ is the singleton, that is, $Q$ is the one-vertex quiver
 and establish a required sink sequence of mutations $\mathcal S_{A_\mathbf N}$ of $A_\infty$.


 Consider the following pair of nested intervals
\[
\{c=r_{m-1}, \cdots , c+d=k_{m-1}+r_{m-1}\}\supset \{a=r_m, \cdots, a+b=k_m+r_m\}.
 \]

We have $c\le a$, $a+b\le c+d$. Hence $k_{m-1}=d>b=k_m$, unless $a=c$ and $b=d$. In the cases
$a=c$ and $b=d$, the claim follows from Proposition \ref{KRcluster}.



The required sequence  
$\mathcal S_{a,b,c,d}$ is a sink green sequence of mutations of the  subquiver $A_{c+d}\subset A_\infty$, and defined as follows.

The sequence $\mathcal S_{a,b,c,d}$
begins with the sequence $\mathcal  S_{A_{a+b-1}}$
\[ 1, 2 ,1, 3, 2, 1\ldots , a+b-1, a+b-2,\cdots, 1.\]
The seed $\Sigma_0 \mathcal  S_{A_{a+b-1}}$ has cluster variables $x^{(a+b)}_{a+b-a'}$, $a'<a+b$ with $g$-vectors
\[
e_{a+b}-e_{ a'}.
\]
Then  we continue by
\begin{equation}\label{seq-ab}
\mathcal S_1:=a+b, a+b-1, \cdots, a+b+1-c,\end{equation}
(a beginning part length $a$ of the sequence $a+b, a+b-1, \cdots, 1$).

Applying this sequence to $\Sigma_0 \mathcal  S_{A_{a+b-1}}$, we get a seed in which
cluster variables whose $g$-vectors have the form
\[
e_{ a+b}-e_{ a'}, \, a<a'<a+b,
\]
do not change (since $a+b+1-c>b=a+b-a$), 
and a new cluster variable at the node $a+ b+1-c$ is equal to
$x^{(a+b+1)}_{a+b+1-c} $ and has $g$-vector
\[
e_{a+b+1}-e_c.
\]
Then we continue with the sequence of mutations 
\[
\mathcal S_2:=a+b+1, a+b, \cdots, a+b+2-c,\]
and get new variables with $g$-vectors 
\[
e_{a+b+2}-e_{c'}, \quad 1\le c'\le c.
\]
Then we apply the sequence 
\[
a+b+2, a+b+1, \cdots, a+b+3-c,\]
and so on
\[\cdots\]
and finally we apply the sequence
\begin{equation}\label{seq-abcd}
c+d-1=a+b+(c+d-a-b), c+d-2, \cdots, d=d+c-c,\end{equation}
and we get new cluster variables with $g$-vectors
\[
e_{c+d}-e_{c'}, \quad 1\le c' \le c. 
\]



From the proof of Proposition \ref{gvectors(n)}, it follows that, after the mutation at the vertex  $a+b+1-c$ in  the sequence (\ref{seq-ab}), we get
the cluster variable $x^{(a+b+1)}_{a+b+1-c}$ whose 
$g$-vector is $e_{a+b+1}-e_c$. 
 
The same reasoning shows that, for $f\in [2, d+c-(a+b)] $, we get
the cluster variable $x^{(a+b+f+1)}_{a+b+f-c}$ after the last mutation following
$a+b+f, \cdots ,a+b+f-c$, and the cluster variables at the vertices labeled
by vertices at $i\le c$ 
do not change.

Finally, we get the cluster variable $x^{(c+d)}_{d}$ with $g$-vector $e_{c+d}-e_c$ after the last mutation
of the sequence $\mathcal S_{a,b,c,d}$.

Thus, for $\widehat{sl_2}$, we get a sink sequence of mutations, such that the quiver $\Sigma_0\mathcal S_{a,b,c,d}$ contains cluster variables with $g$-vectors $e_{a+b}-e_a$ and $e_{c+d}-e_c$. Hence, after change of variables (\ref{tx}) corresponding to these cluster variables are $\widehat{sl_2}$ $q$-characters for $W_{k_{m-1},r_{m-1}+k_{m-1}}$ and $W_{k_{m},r_{m}+k_{m}}$. This finishes the proof for the singleton quiver $Q$.

Since the sequence $\mathcal  S_{a, b, c, d}$ is a sink sequence, for an $ADE$ quiver $Q$, by Proposition \ref{slice1}, for the sequence $Sc(Q)\times \mathcal  S_{a, b, c, d}$ and the mutation at vertex $(v_i,t)=(p_j,t)$ ($(v_i,t)=(q_j,t)$) of the part of the sequence $a+b+f, a+b+f-1, \cdots, a+b+f-c$ of $\mathcal S_{a,b,c,d}$, the neighboring vertex, the head of the incoming edge $((v_i, t-1), (v_i, t))$ is
labeled by the cluster variable $x^{(a+b+f)}_{(v_i,t-1)}$ and other neighboring vertices tails of outgoing edges of the form $(\cdot, t)$
is labeled by $x^{(a+b+f)}_{(\cdot ,t)}$ for $v_i=p_j$  
by $x^{(a+b+f+1)}_{(\cdot ,t)}$ for $v_i=q_j$. The head of the incoming edge
$((v_i, t+1), (v_i, t))$ is labeled by $x^{(a+b+f+1)}_{(v_i ,t+1)}$.
Note that by a similar argument to the proof of Lemma \ref{levelsub},
we see that the sequence has the level property.

 Thus, following the same reasoning as for $\mathcal S_{a,b,c,d}$, we get
 the seed $\Sigma_0 (Sc(Q)\times \mathcal{S}_{a,b,c,d})$ which contains cluster variables $x^{(a+b)}_{v_i,b}$ and $x^{(c+d)}_{v_i,d}$, whose $g$-vectors are
 $e_{(v_i, a+b)}-e_{(v_i, a)}$ and $e_{(v_i, c+d)}-e_{(v_i, c)}$, $v_i\in V_Q$, respectively. Hence, after change of variables, we receive
either the complete or truncated $q$-characters for $W_{k_m,r_m+k_m}^{(v_i)}$ and $W_{k_{m-1},r_{m-1}+k_{m-1}}^{(v_i)}$ depending on $r_m>h/2$ or not.

For the next module $W_{k_{m-2},r_{m-2}+k_{m-2}}^{(v_j)}$ 
 we continue the sequence $\mathcal S_{a,b,c,d}$ by mutation at the following vertices
\[
c+d, c+d-1, \cdots, d+c-r_{m-2}+1,\]
($d+c-r_{m-2}\ge d=d+c-c$, $c=r_{m-1}$)

 \[
 c+d+1, \cdots, d+c+2-r_{m-2},\]
\[ \cdots\]
\[
k_{m-2}+r_{m-2}-1=c+d+(k_{m-2}+r_{m-2}-c-d)-1  , \cdots, \]
\begin{equation}\label{seq-3}
k_{m-2}=k_{m-2}+r_{m-2}-r_{m-2},\end{equation}

Then after the mutation at the vertex $(v_j,k_{m-2})$ in (\ref{seq-3}) we get a cluster variables with $g$-vectors 
\[
e_{(v_j, k_{m-2}+r_{m-2})}- e_{(v_j,r_{m-2})}.
\]
Hence, by Corollary \ref{gq}
we get the statement.

For the next module $W_{k_{m-3},r_{m-3}+k_{m-3}}^{(v_{{j'}})}$, we continue analogously the sequence (\ref{seq-3}) and so on.


  

\hfill $\Box$

\begin{cor}
In the setting of Theorem \ref{KRseed} and its proof, when $Q$ is
then singleton,
$c$-vectors of cluster variables in $\Sigma'$ obtained by applying
mutations after $\mathcal{S}_{A_{a+b-1}}$ in the sequence $\mathcal{S}_{\mathbf{N}}$ form a nested collection, that is, the supports of any pair of $c$-vectors are nested.
\end{cor}
{\em Proof}. From Proposition \ref{gvectors(n)} and the level property, we get that after applying the beginning part (\ref{seq-ab}), $\mathcal S_{A_{a+b-1}}$ the cluster variables at $i=1, \cdots, b$ have $c$-vectors with supports of the form $[a+b-1-i,a+b-1]$. These supports are nested and included in support of other $c$-vectors of this sequence.
Cluster variables obtained by mutations of the remaining part of $\mathcal S_\mathbf N$ gives us a nested sequence with the property
\begin{equation}\label{tight1}
k_i+r_i=k_{i-1}+r_{i-1}-1\mbox{  or   }r_i=r_{i-1}+1.
\end{equation}  
This implies the claim.
\hfill $\Box$  

\begin{defn}
    Two nested collections $\mathbf{N}$ and $\mathbf{N}'$ are in the {\em general position} if it holds $\mathbf{N}\cup \mathbf{N}'$
is nested or \[
\min_{a', [a',b']\in \mathbf N'}a'-\max_{b, [a,b]\in \mathbf N }b >h/2,
\]
or
\[
\min_{a, [a,b]\in \mathbf N}a-\max_{b', [a',b']\in \mathbf N '}b'>h/2.
\]
\end{defn}
\begin{rem} For $\widehat{\mathfrak{sl}_2}$, this condition reads as the union of $q$-strings is not a $q$-string. 

For ADE types, the choice of the difference $h/2+1$ between $\mathbf N$ and $\mathbf N'$ is because in this case the q characters of the modules of $\mathbf N$ and of the modules of $\mathbf N'$ do not share common variables.   
\end{rem}
\begin{thm}\label{KRstring}
    For nested collections $\mathbf{N}$ and $\mathbf{N}'$ in general position there exists a seed $\Sigma_{\mathbf N,\mathbf N'}$ which contains cluster variables with $g$-vectors
    $e_{(v_i, a+b)}-e_{(v_i, a)}$, $[a, a+b]\in \mathbf N$ and
    $e_{(v_i, a'+b')}-e_{(v_i, a')}$, $[a', a'+b']\in \mathbf N'$, $v_i\in Q$.
\end{thm}
From this theorem and Theorem \ref{main} the tensor product of the KR-modules which have $q$-characters for such cluster variables is simple.

{\em Proof}. Suppose we have 
\[ \min_{a', [a', b']\in \mathbf N'}a'-\max_{b, [a, b]\in \mathbf N }b >h/2.
\]
Denote by $a_{\mathbf N'}:=\min_{a, [a,b]\in \mathbf N'}a$.

Then to obtain the required seed, we first apply the sequence of mutations $\mathcal S_\mathbf N$ from Theorem \ref{KRseed}, and let $L$ be such that the sequence $\mathcal S_\mathbf N$ has variables in the set $\{(v, i), v\in Q, i\le L\}$.

Secondly, we apply the sequence $\mathcal S_{\mathbf N'}^{+L}$, which is obtained from $\mathcal S_{\mathbf N'}$ by shifting by $L$ to the right similarly to that $\mathcal S_{A_{[k,K]}}$ is obtained by shifting $\mathcal S_{A_{K-k}}$ to the right by $k$. Then by Proposition \ref{xsub} we get the claim.  \hfill $\Box$

\subsection{Non-nested collections}
We consider the labeling of vertices of a Dynkin quiver of one of $ADE$ types as in Picture 2.

Consider, for $m\ge h/2+1$, and $l_i\ge h/2$, the following collections of intervals \[\mathbf C_{A_n}(l_1, \ldots, l_n;m)=\{[l_i, m+\lfloor\frac i 2\rfloor+1] | i=1, \ldots, n\}\]
 \begin{eqnarray*}
 \mathbf C_{D_n}(l_1, \ldots, l_n;m)&=&\{[l_i, m+\lfloor\frac i 2\rfloor+1]| i=1, \ldots, n-2\}\\
&\cup &\{[l_{n-1}, m+n/2-1], [l_n, m+n/2-1] \}
 \end{eqnarray*}
for even $n$, and 
  \begin{eqnarray*}
 \mathbf C_{D_n}(l_1, \ldots, l_n;m)&=&\{[l_i, m+\lfloor\frac i 2\rfloor+1]| i=1, \ldots, n-2\}\\
 &\cup&[l_{n-1}, m+(n-1)/2], [l_n, m+(n-1)/2]  \}
\end{eqnarray*}
for odd $n$;
\[
\mathbf C_{E_6}(l_1, \ldots, l_6;m)=\{[l_1, m], [l_3, m+1], [l_2, m+1], [l_4, m+1], [l_5, m+2], [l_6, m+2]\},
\]
\[
\mathbf C_{E_7}(l_1, \ldots, l_7;m)=\{[l_1, m], [l_3, m+1], [l_2, m+1], [l_4, m+1], [l_5, m+2], \]
\[[l_6, m+2], [l_7, m+3]\},
\]
\[
\mathbf C_{E_8}(l_1, \ldots, l_8;m)=\{[l_1, m], [l_3, m+1], [l_2, m+1], [l_4, m+1], [l_5, m+2], \]
\[[l_6, m+2], [l_7, m+3], [l_8, m+3]\}.
\]

For each of such a collection, consider the collection of KR-modules whose supports of dominant monomials in $\mathbf C_\cdot(l_1, \cdots, l_n;m)$. For example, for type $A$, the corresponding collection is
\[M_{\mathbf C_{A_n}  (l_1, \cdots, l_n;m)}=\{W^{(v_i)}_{m+\lfloor i/2\rfloor-l_i+2, m+\lfloor i/2\rfloor+1} | i=1,2,\cdots,n \}.
\]

Note that for each of $ADE$ types and for any $m$ and  $l_1<l_2<\cdots<l_n$  the collection $\mathbf C_{\cdot} (l_1, \cdots,l_n;m)$ is not nested.
\begin{thm}\label{nonested}
For each of the types $ADE$ and the corresponding collections $\mathbf C_\cdot (l_1, \cdots,l_n)$, such that all the intervals of the collection are not empty,
there exists a green sequence of mutations $\mathcal S_\mathbf C$ such that the $q$-characters of the KR modules of $M_{\mathbf C_\cdot  (l_1, \cdots,l_n )}$ correspond to some cluster variables of
$\Sigma_0 \mathcal S_\mathbf C$. 
\end{thm}
{\em Proof}. We construct the sequence $\mathcal S_C$ as follows.

For types $A$ and $D$, we begin with the maximal green sequence $Sc(Q)\times \mathcal S_{A_{m-1}}$, then we continue with the following sequence of mutations.

\[
Sc(Q_{\ge 2})\times  (\mathcal S_{A_{m}}\setminus \mathcal S_{A_{m-1}}),\, Sc(Q_{\ge 4})\times  (\mathcal S_{A_{m+1}}\setminus \mathcal S_{A_{m}}),\cdots
\]
Where $Q_{\ge 2l}$ denotes the subquiver of $Q$ with vertices labeled $2l, 2l+1, \cdots$.

For $Q=A_5$, we depicted quivers $Sc(Q)\times \mathcal S_{A_{m-1}}$ (one has to imagine $Sc(Q)\times \mathcal S_{A_{m-1}}$ to be located to the left of the first column in Picture 5 which is $(m-1)$th column with vertices $(v_i, m-1)$, $v_i\in Q$),
$Sc(Q_{\ge 2})\times  (\mathcal S_{A_{m}}\setminus \mathcal S_{A_{m-1}})$,
and $Sc(Q_{\ge 4})\times  (\mathcal S_{A_{m+1}}\setminus \mathcal S_{A_{m}})$ in the pictures below.

\begin{figure}[H]
\centering
  \includegraphics[width=8cm, height=6cm]
  {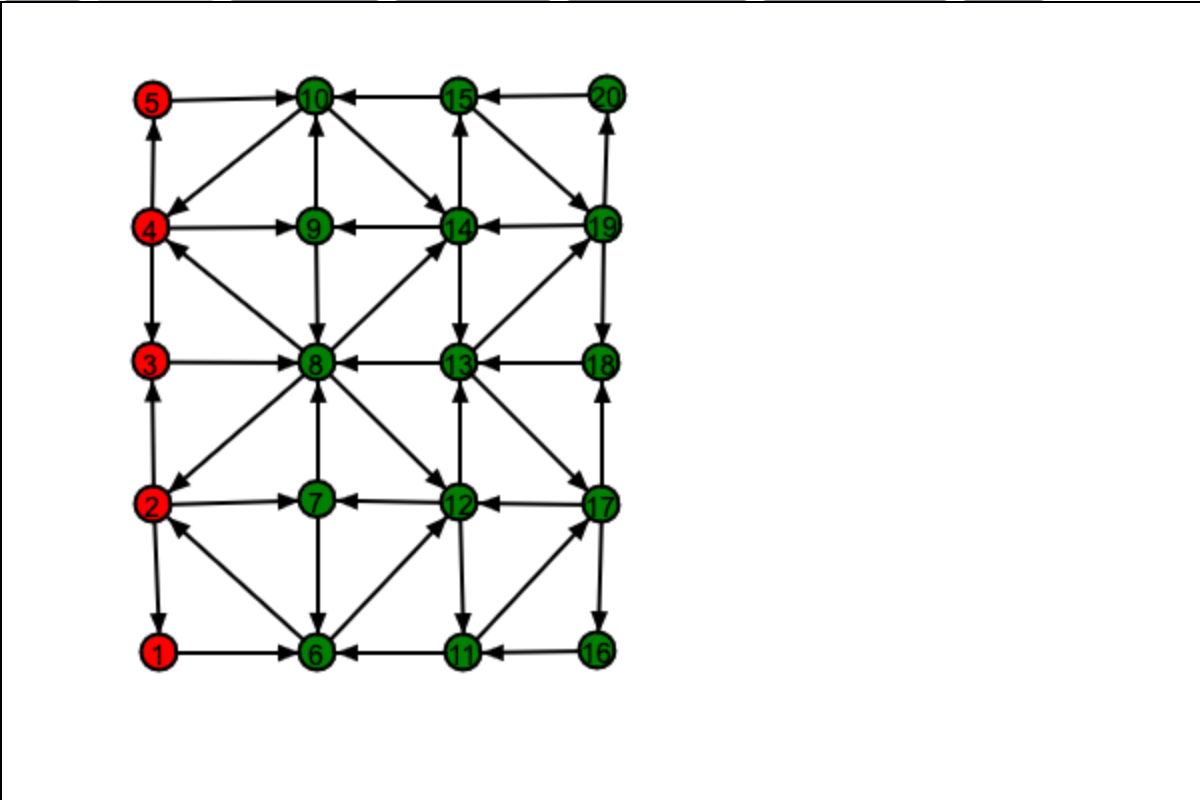}
  \end{figure}  
 $\mbox{ Picture 5\label{nnest1}}$

 \bigskip

 \begin{figure}[H]
\centering
  \includegraphics[width=8cm, height=6cm]
  {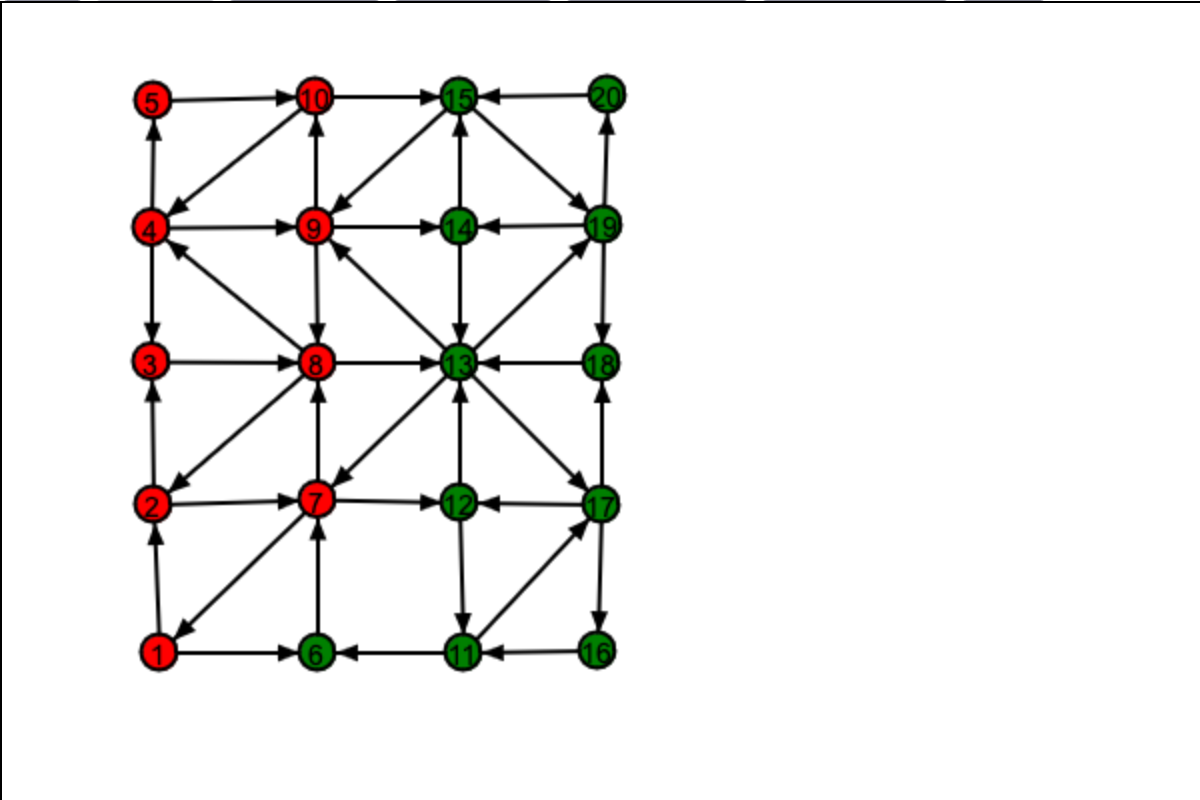}
  \end{figure}  
 $\mbox{ Picture 6\label{nnest1}}$

 \bigskip

 \begin{figure}[H]
\centering
  \includegraphics[width=8cm, height=6cm]
  {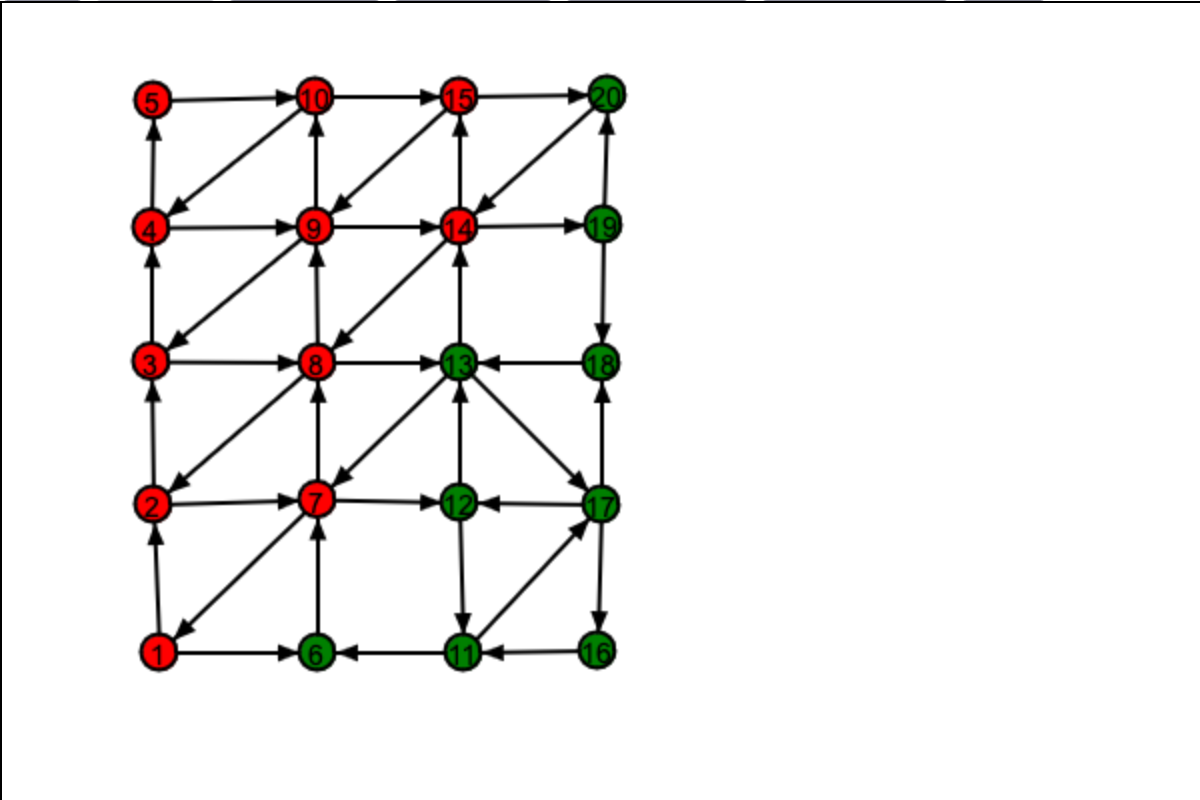}
  \end{figure}  
 $\mbox{ Picture 7\label{nnest1}}$

 \bigskip
For  $Q=E_6$, we expand $Sc(Q)\times \mathcal S_{A_{m-1}}$ by
\[
Sc(Q^{-\{1\}})\times  (\mathcal S_{A_{m}}\setminus \mathcal S_{A_{m-1}}),\, Sc(Q^{-\{1, 3, 4\}})\times  (\mathcal S_{A_{m+1}}\setminus \mathcal S_{A_{m}}),
\]
where $Q^{-L}$ denotes the subquiver of $Q$ without the vertices of the set $L$.

For $Q=E_7$, $Q=E_8$, we expand $Sc(Q)\times \mathcal S_{A_{m-1}}$ by
\[
Sc(Q^{-\{1\}})\times  (\mathcal S_{A_{m}}\setminus \mathcal S_{A_{m-1}}),\, Sc(Q^{-\{1, 3, 4\}})\times  (\mathcal S_{A_{m+1}}\setminus \mathcal S_{A_{m}}), \,\]
\[Sc(Q^{-\{1, 2, 3, 4, 5\}})\times  (\mathcal S_{A_{m+2}}\setminus \mathcal S_{A_{m+1}}).
\]
For type $A$, the set of $g$-vectors of cluster variables of the seed $\Sigma_0 \mathcal S_\mathbf C$ contains the set $\{e_{(v_i,m+\lfloor i/2\rfloor+1)}-e_{(v_i,l_i)}|
[l_i,m+\lfloor i/2\rfloor+1]\in \mathbf C(l_1, \cdots , l_n, m)\}$, and corresponding to $\mathbf C(l_1, \cdots,l_n;m )$ sets of $g$-vectors for other types.

This follows from proving the level property, namely the mutation at the source vertex $(p_i, k)$ or $(q_j, k)$ is a mutation at a vertex with
two incoming edges with head vertices $(p_i, k\pm 1)$ or $(q_j, k\pm 1)$, and the cluster variable with the vertices $(p_i, k-1)$ or $(q_j, k- 1)$
is computed in the previous step, and the cluster variable with the vertices $(p_i, k+1)$ or $(q_j, k+ 1)$
are computed on the current step. For cluster variables at the tails of the outgoing edges
of the form $(\cdot, k)$ are computed in the previous step for $p_i$'s and are computed in the current step for $q_j$'s. Since in each step we start from a source vertex $p_i$, the current step computed cluster variables at these vertices uses the previous step cluster variables at the tails of outgoing edges $(\cdot, k)$. Since $q_j$'s from the sequence are located 'above' the first $p_j$ in $Sc(Q_{\ge 2\cdot})$
we see that the cluster variables which we need for computation of a new variable are yet computed at the current step. 
 
 Since $v_2=2$ is a source vertex, for computing the cluster variables in $(\Sigma_0 (Sc(Q)\times \mathcal S_{A_{m-1})})Sc(Q_{\ge 2})\times (\mathcal S_{A_m}\setminus \mathcal S_{A_{m-1})}$ we get the same cluster variables labeled $(v_i, j)$, $i\ge 2$, $j\le m$ as in the seed $\Sigma_0 (Sc(Q)\times \mathcal S_{A_{m}})$, since we do not need their computation to recalculate the values in the cluster variables labeled $(v_1, j)$ in $\Sigma_0 (Sc(Q)\times \mathcal S_{A_{m-1}})$.

 Since $4$ is the next source vertex, the same property holds for $Q_{\ge 4}$. Namely, for the seed \[(\Sigma_0 (Sc(Q)\times \mathcal S_{A_{m-1})})Sc(Q_{\ge 2})\times (\mathcal S_{A_m}\setminus \mathcal S_{A_{m-1}})\times 
 Sc(Q_{\ge 4})\times (\mathcal S_{A_{m+1}}\setminus \mathcal S_{A_{m})}\] we get the same cluster variables labeled $(v_i, j)$, $i\ge 4$, $j\le m+1$ as in the seed $\Sigma_0 (Sc(Q)\times \mathcal S_{A_{m+1}})$, since we do not need their computation to recalculate the values in the cluster variables labeled $(v_1, j)$, $(v_2, j)$, and $(v_3, j)$ in $\Sigma_0 (Sc(Q)\times \mathcal S_{A_{m}})$
 and so on.
 
For other types, the same property of mutations at source vertices gives us the same property of cluster variables in the seeds of parts of $\mathcal S_\mathbf C$.

Because of this, we obtain $g$-vectors of the resulting seed containing the sets corresponding to $\mathbf C(l_1, \cdots,l_n; m)$ by Proposition
\ref{gvectors(n)}.

\hfill $\Box$

\begin{con}
The tensor product
\[\otimes _{W\in \mathcal M_\mathbf C}W
\]
is simple.
\end{con}
\begin{rem}
 The intervals $[l,m+i]$ and $[l', m+i+1]$  correspond to a pair of commuting 
 KR-modules, while, for $l<l'$, the intervals $[l,m+i]$ and $[l', m+i+2]$  do not correspond to a pair of commuting  
 KR-modules. Therefore, non-nested collections of the conjecture do not correspond to chains of $i$-boxes. 
\end{rem}


\section{Cluster Donaldson-Thomas transformation}\label{DT}
 \subsection{Definition and motivation example}
 
 Let $R$ be a quiver without frozen and $\Sigma=(\mathbf x, R)$ be an initial seed of  a cluster algebra $\mathcal R$.
 Recall that if the quiver $R$ possesses a maximal green (reddening) sequence $\mathcal S$ then there is the unique isomorphism $\sigma:\underline R\cong \overline R{\mathcal S}$ between the co-framed quiver $\underline R$ and the quiver $\overline R{\mathcal S}$ fixing the frozen vertices and sending a non-frozen
 vertex $v$ to $\sigma (v)$ (see Section \ref{basicfram}). 

 The Donaldson-Thomas cluster transformation is a birational map that sends the initial cluster variables to the cluster variables obtained by applying a maximal green (reddening) sequence composed with the isomorphism $\sigma$.
 
 Namely, a birational map
 \begin{equation}\label{DTdef}
 DT_\Sigma: \mathbb C^{V_R}\to \mathbb C^{V_R},\quad  
 \mathbf {x}\to ((\mu_\mathcal S \mathbf{x})_{\sigma v})_{v\in V_R}.
     \end{equation}
 is called the {\em cluster Donaldson-Thomas transformation } in \cite{GonS} (see also \cite{Wang}). The cluster Donaldson-Thomas transformation is an automorphism of the cluster algebra $\mathcal R$.

 Note that if $R$ is an ice quiver, we define the cluster DT-transformation as a birational map sending cluster variables labeled by mutable vertices by (\ref{DTdef}) and being identical at frozen vertices.
  
 The cluster DT-transformation is invariant of the cluster algebra in the following sense. Firstly, $DT_\Sigma$ does not depend on a reddening sequence for $\overline R$. Secondly, let, for $\Sigma=(\mathbf x, R)$ and a mutable vertex $v\in V^{mt}_R$, $\Sigma'=\mu_v\Sigma$ be a seed obtained by mutation in $v$, and $\mathcal S$ be a reddening sequence for $R$. Then (see \cite{GonS})\, 
 the sequence $\mathcal S'=(v, \mathcal S, \sigma (v))$ is a reddening sequence for $R'$ and there holds 
 \begin{equation}
  DT_{\Sigma'}=\mu_{\sigma (v)}\cdot DT_\Sigma \cdot \mu_v,
     \end{equation}
where $\mu_v$ and $\mu_{\sigma (v)}$ are considered as birational maps.

 \begin{ex} Consider a quiver $R=A_{n+1}$ and
 a seed $\Sigma=(\mathbf{x},R)$. That is, we consider the ice quiver $A_{n+1}$ with edge orientation $(i+1, i)$ $i=1, \ldots, n$, the frozen vertex $n+1$, and the maximal green sequence $\mathcal S_{A_n}$. Then by Proposition \ref{gvectors(n)}, and since the automorphism $\sigma$ is the reverse permutation $\sigma (a)=n-a+1$, the cluster DT-transformation for the cluster algebra with the ice quiver $A_{n+1}$ takes the form
  \begin{equation}\label{DT-A}
  DT_{\Sigma}(\mathbf{x})_a=x_{n+1-a}^{(n+1)}=\frac{x_{n+1}}{x_{a}}(1+ y_a+y_ay_{a+1}+\cdots +y_ay_{a+1}\cdots y_n)
\end{equation}
where \[
y_j=\frac{x_{j-1}}{x_{j+1}}, \, j=1, \ldots, n.
\]

Because of (\ref{Fform2}), we can regard (\ref{DT-A}) as the $q$-character
\begin{equation}\label{qDTA}
DT_{\Sigma}(\mathbf{x})_a=Y_{q^{n+1}}\cdots Y_{q^{a+1}}(1+ A_a^{-1}+A_a^{-1}A_{a+1}^{-1}+\cdots +A_a^{-1}A_{a+1}^{-1}\cdots A_n^{-1}),\end{equation} where \[
Y_{q^k}:=\frac{x_k}{x_{k-1}}, 
\,\, A_k^{-1}=\frac 1{Y_{q^k}Y_{q^{k+1}}}. 
\] To obtain the {\em DT}-cluster transformation for $A_n$ without freezing, we have to set $x_{n+1}=1$ in (\ref{DT-A}).  
\end{ex}

The cluster DT-transformation of the form (\ref{qDTA}) is a map sending the cluster variables of the initial seed to the $q$-characters of the KR modules.

For triangular products of Dynkin quivers, one of which is of A-type,  Proposition \ref{KRcluster}\, leads to a relation between the cluster Donaldson-Thomas transformation and $q$-characters (truncated) of Kirillov-Reshetikhin modules. 

 \subsection{DT-transformation for triangular products
 and 
 algorithms computing $q$-characters
 }

 We consider as above Dynkin quiver $Q$ of $ADE$-type with alternating orientation of edges for simply-laced cases 
 and the triangular product $Q\boxtimes A_{m+1}$, 
 where $A_{m+1}$ is the iced $A_{m+1}$ quiver with the frozen vertex $m+1$. Therefore, $Q\boxtimes A_{m+1}$ has frozen vertices of the form
 $(v_i,m+1)$, $i\in [n]$.

 \begin{thm}\label{DT0} 
For the cluster algebra $\mathcal A_{Q;m+1}$ with the initial seed $\Sigma =(\mathbf x, Q\boxtimes A_{m+1})$ with frozen vertices $(v_i, m+1)$, $v_i\in V$, 
 the cluster DT-transformation takes the form
 \begin{equation}\label{DT-0}
x_{v_i, j}\to \frac{x_{(v_i, m+1)}}{x_{(v_i, m+1-j)}} F_{x^{(m+1)}_{(v_i, m+1-j)}}(\mathbf  y)
\quad j=1, \ldots, m
\end{equation}

For the initial seed $\Sigma =(\mathbf x, Q\boxtimes A_{m})$ without frozen 
the cluster DT-transformation takes the form
 \begin{equation}\label{DTunfrozen}
x_{v_i, j}\to \frac{1}{x_{(v_i, m+1-j)}} F_{x^{(m+1)}_{(v_i, m+1-j)}}(\mathbf  y)
\quad j=1, \ldots, m
\end{equation}
After change of coordinates (\ref{tx}) in (\ref{DT-0}), we get the cluster DT-transformation is of the form 
 \begin{equation}\label{DT1}
x_{v, j}\to 
\chi_q (W^{(v)}_{m+1-j,m+1}),
\quad j=1, \ldots, m
\end{equation}
where $\chi_q (W^{(v)}_{m+1-j,m+1})$ the $q$-character of a KR-module, and  for $j\ge 1+ h/2$, $ \chi_q$ denotes the $q$-character of $W^{(v)}_{m+1-j,m+1}$, and, for  $j\le h/2$, the truncated $q$-character (due to Remark \ref{trunc1})
\end{thm}

{\em Proof}.  (\ref{DT-0}) follows from (\ref{sep-product1}) in Proposition \ref{slice1}.
Because frozen are of the form $(v_i, m+1)$, $i\in [n]$, then by freezing, we get that (\ref{DT-0}) turns into (\ref{DTunfrozen}). 

Since the automorphism $\sigma$ is the reverse permutation on $A_m$, (\ref{DT1}) follows from Proposition \ref{KRcluster}.
\hfill $\Box$

For triangular products,  the relation (\ref{DT1}) between the cluster DT-transfor\-mation and $q$-characters of Kirillov-Reshetikhin modules allows
us to use
the Frenkel-Mukhin algorithm or Nakajima algorithm, which compute the $q$-characters of KR-modules, for computing the cluster DT-transformation. (It is easy to modify these algorithms for computing truncated $q$-characters as well.) 

Note that this algorithm is much faster than the original cluster computations of DT-transformation, because the cluster transformation requires division of polynomials in many variables.

 \subsection
{Cluster DT-transformation for  double Bruhat cells
}
There is a structure of cluster algebra on the coordinate rings of double Bruhar cells.

For $u, v\in W$, the double Bruhat cell $B^{u,v}$ is defined as
\begin{equation}\label{double1}
B^{u,v}=B\overline u B\cap B^-\overline v B^-
\end{equation}

We consider varieties 
 $B^-_{w_0}:=B^-\cap U\overline{w_0}U $ and \[
U_{w_0}=U\cap B^-\overline w_0 B^-\]

Let $Q$ be a Dynkin quiver for the Cartan matrix $A$. Then the alternating orientation of edges defines a specific reduced decomposition of $w_0$, the power of the bipartite Coxeter element. 


Recall that $B^-_{w_0}=B^-\cap Uw_0U$.
 By specifying Theorem \ref{DT0} to the case $m+1=h/2$, we get a case of the cluster algebra for $\mathbb C[ B^-_{w_0}]$. To avoid technicalities, we consider the type $A$ with odd number $n$ of vertices in order have an even $h=n+1$. 
 
 In fact, in  such a case of type $A$ and other simply-laced types $D$, $E$, we regard the triangular product 
 \[
 Q\boxtimes A_{h/2}
 \]
 as an initial ice quiver for $\mathbb C[ B^-_{w_0}]$ corresponding to a reduced decomposition
 \[
 w_0=(C^{op})^{h/2},
 \]
where $C= \prod_jq_j\prod_ip_i$ is a bipartite Coxeter word and $C^{op}= \prod_ip_i\prod_jq_j$.

The frozen vertices take the form $(p_i, {h/2})$ or $(q_j, {h/2})$, $i=1, \ldots, k$, $j=1, \dots, l$.

 Then by Theorem \ref{DT0}, for the initial seed $\Sigma_0=(\mathbf x, Q\boxtimes A_{h/2})$ of the cluster algebra for $\mathbb C[ B^-_{w_0}]$, the cluster DT-transformation is defined by specification of 
 (\ref{DT-0}) to $m+1=h/2$. The $q$-characters in (\ref{DT1}) become all truncated.
 
 \section*{Acknowledgments}
 The authors thanks Masaki Kashiwara, Bernard Leclerc and Hironori Oya for useful discussions. G.K. thanks Sophia University for hospitality during visit in July 2024. Y.K. is
 supported by JSPS KAKENHI Grant Number JP24K22825.
 The research of T.N. is supported by
Grant-in-Aid for Scientific Research (C) 20K03564, 
Japan Society for the Promotion of Science.


\begin{thebibliography}{20}
  
 \bibitem{BZ99} A.Berenstein and
A.Zelevinsky, Tensor product multiplicities, canonical bases and totally positive varieties, 
Invent. Math. 143, no. 1, 77--128 (2001).

 
 \bibitem{BFZ3} A.Berenstein, S.Fomin, and A.Zelevinsky, Cluster algebras III: Upper bounds and double Bruhat cells, Duke Math. J. 126 (2005), no. 1, 1–52.
 

 \bibitem{BDP} T. Br\"ustle, G. Dupont, and M. P\'erotin. “On maximal green sequences”. International Math- ematics Research Notices 2014.16 (2014), pp. 4547–4586.
 \bibitem{CKQ} P. Cao, F. Qin, and B. Keller,  The valuation pairing on an upper cluster algebra, 2022,
 https://webusers.imj-prg.fr/~bernhard.keller/publ/index.html
\bibitem{CP}
Vyjayanthi Chari and Andrew Pressley. Minimal affinizations of representations of quantum groups: the
simply laced case. J. Algebra, 184(1), 1-30, 1996.

\bibitem{CP1} V. Chari , A. Pressley, Quantum affine algebras, Comm. Math. Phys. 142 (1991), 261–283.

\bibitem{FHOO} Ryo Fujita, David Hernandez, Se-jin Oh, and Hironori Oya, Isomorphisms among quantum Grothendieck rings and propagation of positivity, arXiv:2101.07489 

  \bibitem{Keller} L. Demonet and B. Keller. “A survey on maximal green sequences” (2019). arXiv: 1904.09247.

 \bibitem{DWZ} H. Derksen, J. Weyman, and A. Zelevinsky, Quivers with potentials and their representations II: applications to cluster algebras,  Journal of the American Mathematical Society, 23/3 (2010), pp. 749--790.

 \bibitem{FZ4} S. Fomin, and A. Zelevinsky, Cluster algebras IV: Coefficients, Comp. Math. 143 (2007), 112--164.

 \bibitem{FM} E.Frenkel and E.Mukhin, Combinatorics of q-characters of finite-dimensional repre- sentations of quantum affine algebras, Comm. Math. Phys. 216 (2001), 23–57.
 
\bibitem{FR} E.Frenkel and N. Reshetikhin,  The $q$-characters of representations of quantum affine algebras, Recent
developments in quantum affine algebras and related topics, {\em Contemp. Math.} 248 (1999), 163--205.


\bibitem{Tsys} R. Inoue, O. Iyama, A. Kuniba, T. Nakanishi and J. Suzuki, Periodicities of T-systems
and Y-systems, Nagoya Math. J.197(2010), 59--174.

 \bibitem{GK} V.Genz and G.Koshevoy,  Maximal green sequences for triangle products,  {\em S\'eminaire Lotharingien de Combinatoire 
 FPSAC 2021}, Issue 86B
 
 \bibitem{GKSadv}  V.Genz, G. Koshevoy, and B.Schumann 
 Polyhedral parametrizations of canonical bases \& cluster duality, {\em Advances in Mathematics}, 2020, 369, 107178 
 
 \bibitem{GonS} A. Goncharov and L. Shen. Donaldson-Thomas transformations of moduli spaces of G-local systems, Advances in Math. 327: 225--348 
 
 \bibitem{Kellerlin} Giovanni Cerulli Irelli, Bernhard Keller, Daniel Labardini-Fragoso and
Pierre-Guy Plamondon, 
Linear independence of cluster monomials for
skew-symmetric cluster algebras, Compositio Math. 149 (2013) 1753–1764

\bibitem{Tsys} R. Inoue, O. Iyama, A. Kuniba, T. Nakanishi ad J. Suzuki, Periodicities of T-systems
and Y-systems, Nagoya Math. J.197(2010), 59--174.

 \bibitem{GK} V.Genz and G.Koshevoy,  Maximal green sequences for triangle products,  {\em S\'eminaire Lotharingien de Combinatoire 
 FPSAC 2021}, Issue 86B
 
 \bibitem{GKSadv}  V.Genz, G. Koshevoy, and B.Schumann 
 Polyhedral parametrizations of canonical bases \& cluster duality, {\em Advances in Mathematics}, 2020, 369, 107178 
 
 \bibitem{GonS} A. Goncharov and L. Shen. Donaldson-Thomas transformations of moduli spaces of G-local systems, Advances in Math. 327: 225--348 
 
 \bibitem{Hernandez} D.Hernandez, Simple tensor products, Invent math (2010) 181: 649 -- 675,  arXiv: 0907.3002

 \bibitem{Hernandez1} D.Hernandez, The Kirillov-Reshetikhin conjecture and solutions to T-systems. {\em J.Reine. Angew. Math.} 596 (206), 63--87
 
 \bibitem{HL} Hernandez, David; Leclerc, Bernard. A cluster algebra approach to $q$-characters of Kirillov-Reshetikhin modules. J. Eur. Math. Soc. (JEMS) 18 (2016), no. 5, 1113--1159.
 arXiv:1303.0744
 

   \bibitem{HL1} Hernandez, David; Leclerc, Bernard. Quantum affine algebras and cluster algebras. Interactions of quantum affine algebras with cluster algebras, current algebras and categorification—in honor of Vyjayanthi Chari on the occasion of her 60th birthday, 37--65, Progr. Math., 337, Birkhäuser/Springer, Cham, [2021], ©2021.  MR4404353
   arXiv:1902.01432v2


 \bibitem{Kac} V.G. Kac, Infinite dimensional Lie algebras, Cambridge University Press, Cambridge, Third edition, 1990 

 
 \bibitem{HL} D. Hernandez and B.Leclerc, A cluster algebra approach to $q$-characters of Kirillov-Reshetikhin modules, 
 J. Eur. Math. Soc. (JEMS) 18 (2016), no. 5, 1113–1159.
 
 
  \bibitem{KKN1} 
 Y. Kanakubo, G. Koshevoy, and T. Nakashima,
 An algorithm for Berenstein-Kazhdan decoration functions and trails for minuscule representations,
Journal of Algebra, 608(2022), pp. 106–142
 
 \bibitem{KKN2} 
 Y. Kanakubo, G. Koshevoy, and T. Nakashima, An algorithm for Berenstein-Kazhdan decoration functions and trails for classical Lie algebras, 
 Int. Math. Res. Not. IMRN (2024), no. 4, 3223–3277.

\bibitem{KKKO}
Kang, S.-J., Kashiwara, M., Kim, M., Oh, S.-j.,
Simplicity of heads and socles of tensor products,
Compos. Math. 151 (2015), no. 2, 377–396.


\bibitem{KashiwaraKorea} Kashiwara, M., Kim, M., Oh, S.-j., Park, E.: Monoidal categorification and quantum affine algebras II, Inventiones mathematicae (2024) 236: 837--924

 \bibitem{Kellerperiod} B.Keller, The periodicity conjecture for pairs of Dynkin diagrams, 
 Annals of Mathematics, volume 177, issue 1 (January 2013). 
 
 
 \bibitem{Mul} Greg M\"uller,  The existence of a maximal green sequence is not invariant under quiver mutation, Electron. J. Combin. 23 (2016), no. 2, Paper 2.47
 
 \bibitem{Nakajima} H. Nakajima, Quiver varieties and t-analogs of q-characters of quantum affine algebras. Annals of mathematics, pages 1057--1097, 2004.
 
 \bibitem{NZ} T.Nakanishi and A.Zelevinsky, On Tropical Dualities in Cluster Algebras, Algebraic groups and quantum groups, 
565 (2012),  217--226
 
 
 \bibitem{OT}  S. Oh and  T. Scrimshaw,  Simplicity of tensor products of Kirillov--Reshetikhin modules: nonexceptional affine and $G$ types,  	arXiv:1910.10347

 
 
 \bibitem{Qin} F. Qin, Triangular bases in quantum cluster algebras and monoidal categorification conjectures, Duke Math. J.166(2017), no. 12, 2337–2442
 
 \bibitem{Xie} D.Xie, BPS spectrum, wall crossing and quantum dilogarithm identity, Adv. Theor. Math. Phys. 20 (2016),
no. 3, 405–524
arXiv:1211.7071
 
 \bibitem{Wang} D. Weng, F-polynomials of Donaldson-Thomas transformations, arXiv:2303.03466
 \end{thebibliography}
\end{document}